# Generalizations of Bricard octahedra

Gerald D. Nelson
nelso229@comcast.net

**Abstract**

We construct five types of polyhedra by generalizing the description of Bricard octahedra and applying the generalizations to polyhedral suspensions. The resulting polyhedra are flexible, are of genus 0, exhibit self-intersections, have zero oriented volume, have dihedral angles all of which are non-constant during flexion and are of indefinite size, the smallest of which are dodecahedra with eight vertexes.

## 1. Introduction

Bricard octahedra are flexible polyhedra in which dihedral angles can change while the geometric characteristics (angles and edge lengths) of their triangular faces remain unchanged. These polyhedra, of which there are three types, were discovered by Bricard in 1897 [1]. The motion of vertexes in Bricard octahedra can be respectively characterized as axial, planar or moving between two co-planar positions. Alternately there are geometric characteristics unique to each type that can be used to describe the octahedra. For the first two types there are six pairs of equal length edges, appropriately arranged; for the third type opposite face angles at all vertexes are either equal or supplementary in value.

In this paper we show that five distinct families of flexible polyhedra can be constructed from generalizations of these geometric characteristics. The generalizations are based upon geometric properties of the various distinct tetrahedral caps (open four facetted polyhedral surfaces) that occur in Bricard octahedra and are formulated to describe polyhedral suspensions (polyhedra having the combinatorial structure of dipyramids) that are comprised of two larger caps each having N, an even integer $\geq 6$, faces as illustrated in Fig. 1. The two apical, or cap, vertexes $\mathbf{u}$ and $\mathbf{w}$, are of index N while the vertexes at the base of the cap, $\mathbf{v_1}$, $\mathbf{v_2}$,…$\mathbf{v_N}$, are all of index 4.

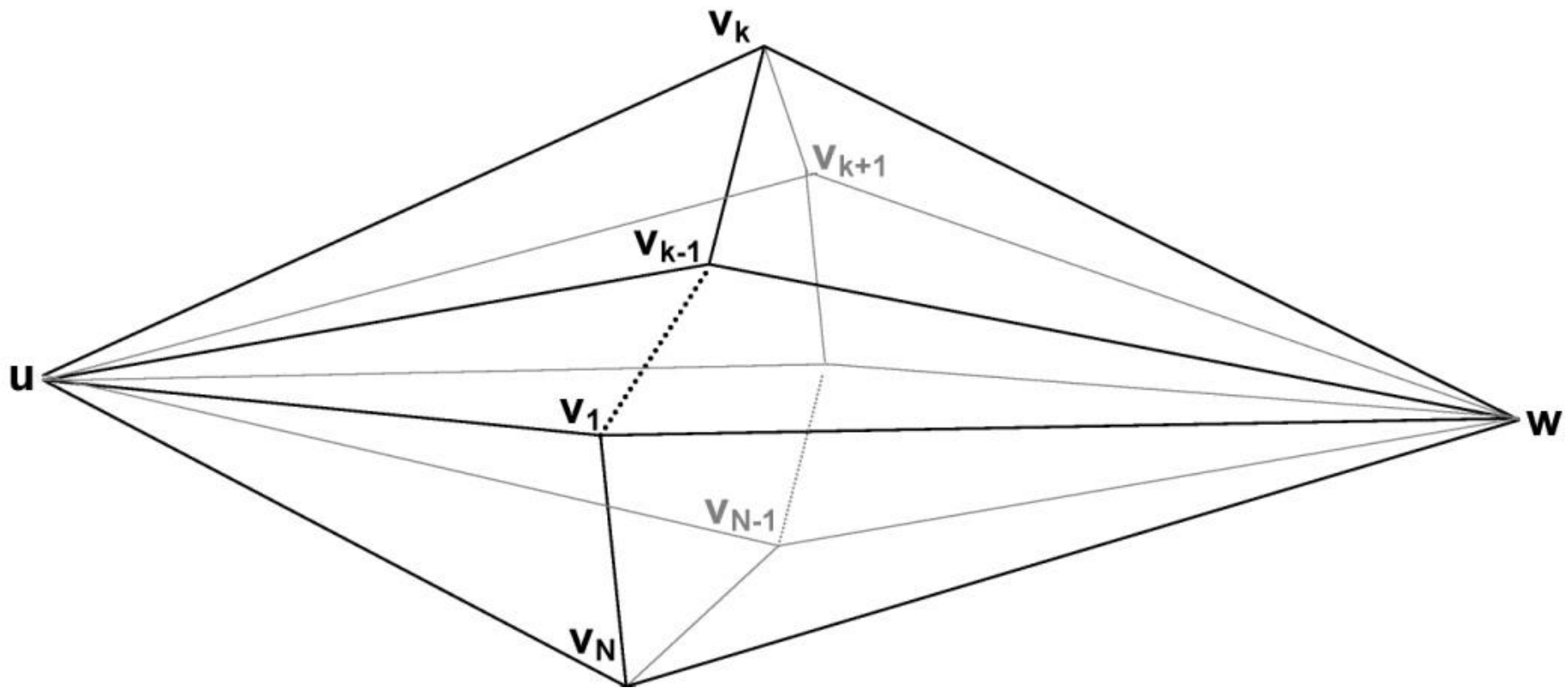

**Figure 1.** Polyhedral Suspension with apical vertexes of Index N

The generalizations lead to flexible polyhedra of genus 0 that have N+2 vertexes and 2N faces. Like the Bricard octahedra from which they are drawn these polyhedra have faces that exhibit self-intersections and under flexion all of their dihedral angles are non-constant.

Polyhedra of this form have been previously and widely studied. In [2] a parameterization of flexible suspensions from the elliptic Wierstrasse function is demonstrated. [3] provides two theorems on the rigidity of suspensions.

Some necessary conditions for flexibility are presented in [4]. In [5] the construction of a new flexible dodecahedron (N=6, Fig. 1) is described and its properties studied. However, explicit geometric descriptions of flexible suspensions of the type provided in this paper have not been previously reported.

We proceed as follows: A classification that is based upon five distinct sub-types of tetrahedral caps that are found in Bricard octahedra is defined (Sec. 2) and some notation developed to describe suspensions (Sec. 3). Generalizations are made for each of the five identified sub-types and applied to suspensions with even values of N>4 leading to the parametric specification of five families of suspensions (Secs. 4-9). An approach to the construction of these suspensions is then described along with the description of a test for flexibility (Sec. 10). A discussion (Sec. 11) of flexible suspensions that have been constructed using the generalizations is included. We conclude with conjectures (Sec. 12) that the suspensions described are indeed flexible for all even N>4. Appendix A describes an approach for defining all parameters of a Bricard octahedron of the third type. Appendix B describes an approach for defining all parameters of a flexible suspension that results from a generalization of a Bricard octahedron of the third type. Appendix C develops equations that describe the deformation of the caps of a suspension. Appendix D contains data tables that parameterize several example flexible suspensions and images of selected examples.

## 2. Bricard Octahedra Caps

Five distinct tetrahedral caps can be found in Bricard octahedra and can be designated as follows: I-OEE, II-AEE, II-OEE, III-OAE and III-OAS. Here the Roman numeral refers to the type of the associated Bricard octahedron and the abbreviated identifiers refer to geometric characteristics as described in Table I and in the following summary descriptions. These five sub-types have been previously described [6].

| Identifier | Definition |
|---|---|
| OEE | Opposite Edges Equal |
| AEE | Adjacent Edges Equal |
| OAE | Opposite Angles Equal |
| OAS | Opposite Angles Supplementary |

**Table I.** Geometric Characteristics

In I-OEE caps opposite edges on the base of the cap are of equal length, for example $|\mathbf{v_1 v_2}|=|\mathbf{v_3 v_4}|$. II-AEE caps have adjacent edge lengths on the base of the cap that are of equal length, for example $|\mathbf{v_1 v_2}|=|\mathbf{v_4 v_1}|$. For both I-OEE and II-AEE, edge lengths at the apical vertex are not constrained by any specific relationships except that they cannot all be of the same length. II-OEE caps have OEE edge lengths on the base of the cap while edge lengths at the apical vertex are constrained so that opposite edges have equal length, for example $|\mathbf{uv_1}|=|\mathbf{uv_3}|$. In III-OAE caps the face angles at the apical vertex have opposite angles that are equal while in III-OAS caps these angles have supplementary values. With the exception of the II-OEE sub-type, the parameters of one cap are sufficient to completely determine the second cap and completely define the flexible octahedron formed from the two caps. Two additional parameters associated with the second apical vertex are required in the exceptional case.

When applied to Bricard octahedra there is considerable ambiguity in the manner in which the above classifications are used and several possibilities exist for labeling of the vertexes. Opposite vertexes of Bricard octahedra are of the same sub-type; thus there always six choices for the two apical vertexes $\mathbf{u}$ and $\mathbf{w}$ for example. For Bricard octahedra of the first type any pair of opposite vertexes may be designated as the apical vertexes and all choices result in the I-OEE sub-type designation. For Bricard octahedra of the second type two pairs of opposite vertexes may be designated as the apical vertexes for the II-AEE sub-type and one pair for the II-OEE sub-type. Similarly, for Bricard octahedra of the third type two pairs of opposite vertexes may be designated as the apical vertexes for the III-OAE sub-type and one pair for the III-OAS sub-type.

It is possible to define Bricard octahedra of the first and second types that also have the properties of the third type; for example, an octahedron having axial symmetry and which also has two flat positions is possible. For the purpose of uniqueness we designate these special cases as sub-types III-OAE or III-OAS as appropriate.

In addition to the edge length and face angle relationships that parametrically characterize Bricard octahedra their dihedral and solid angles have special relationships that are exhibited when the octahedra are flexed. In all three types, edges that are opposite of one another, for example $\mathbf{uv_1}$ and $\mathbf{wv_3}$, have dihedral angles that are either equal or conjugate in value and in the third type, edges that are opposite of one another at a vertex also have this relationship. Solid angles of vertexes that are opposite of one another are "space filling" in the sense that their solid angles sum to $4\pi$ and in the second and third type there are vertexes that exhibit a constant solid angle value = $2\pi$.

## 3. Some Notation

In the following sections it is convenient to discuss characterizations that are based upon edge lengths, face angles, dihedral angles and solid angles; for that purpose we make the following definitions with reference to Fig. 1. N is an even integer that always refers to the number of vertexes at the base of a cap and N=2M for some M>1. Edge lengths at the apical vertex $\mathbf{u}$ are defined by $l_k=|\mathbf{uv_k}|$, at vertex $\mathbf{w}$ by $m_k=|\mathbf{wv_k}|$ for k=1..N and on the base of the caps, $L_N=|\mathbf{v_N v_1}|$ and $L_k=|\mathbf{v_k v_{k+1}}|$ for k=1..N-1. Face angles at the vertex $\mathbf{u}$ are defined by $\alpha_N=\angle\mathbf{v_N u v_1}$ and $\alpha_k=\angle\mathbf{v_k u v_{k+1}}$; at vertex $\mathbf{w}$, $A_N=\angle\mathbf{v_1 w v_N}$ and $A_k=\angle\mathbf{v_{k+1} w v_k}$ and at the non-apical vertexes $\beta_N=\angle\mathbf{v_1 v_N u}$, $\gamma_N=\angle\mathbf{u v_1 v_N}$, $B_N=\angle\mathbf{w v_N v_1}$, $\Gamma_N=\angle\mathbf{v_N v_1 w}$, $\beta_k=\angle\mathbf{v_{k+1} v_k u}$, $\gamma_k=\angle\mathbf{u v_{k+1} v_k}$, $B_k=\angle\mathbf{w v_k v_{k+1}}$ and $\Gamma_k=\angle\mathbf{v_k v_{k+1} w}$ for k=1..N-1. This notation is illustrated in Fig. 2. Dihedral angles (not illustrated) along the edges at vertex $\mathbf{u}$ are defined by $\delta_k=\angle\mathbf{uv_k}$, along the edges at vertex $\mathbf{w}$ are $\Delta_k=\angle\mathbf{wv_k}$ for k=1..N while at the base of the caps, $\varepsilon_N=\angle\mathbf{v_N v_1}$ and $\varepsilon_k=\angle\mathbf{v_k v_{K+1}}$ for k=1..N-1. Solid angles are defined by $\sigma_k$ for vertex $\mathbf{v_k}$ and by $\sigma_u$ and $\sigma_w$ for the respective apical vertexes.

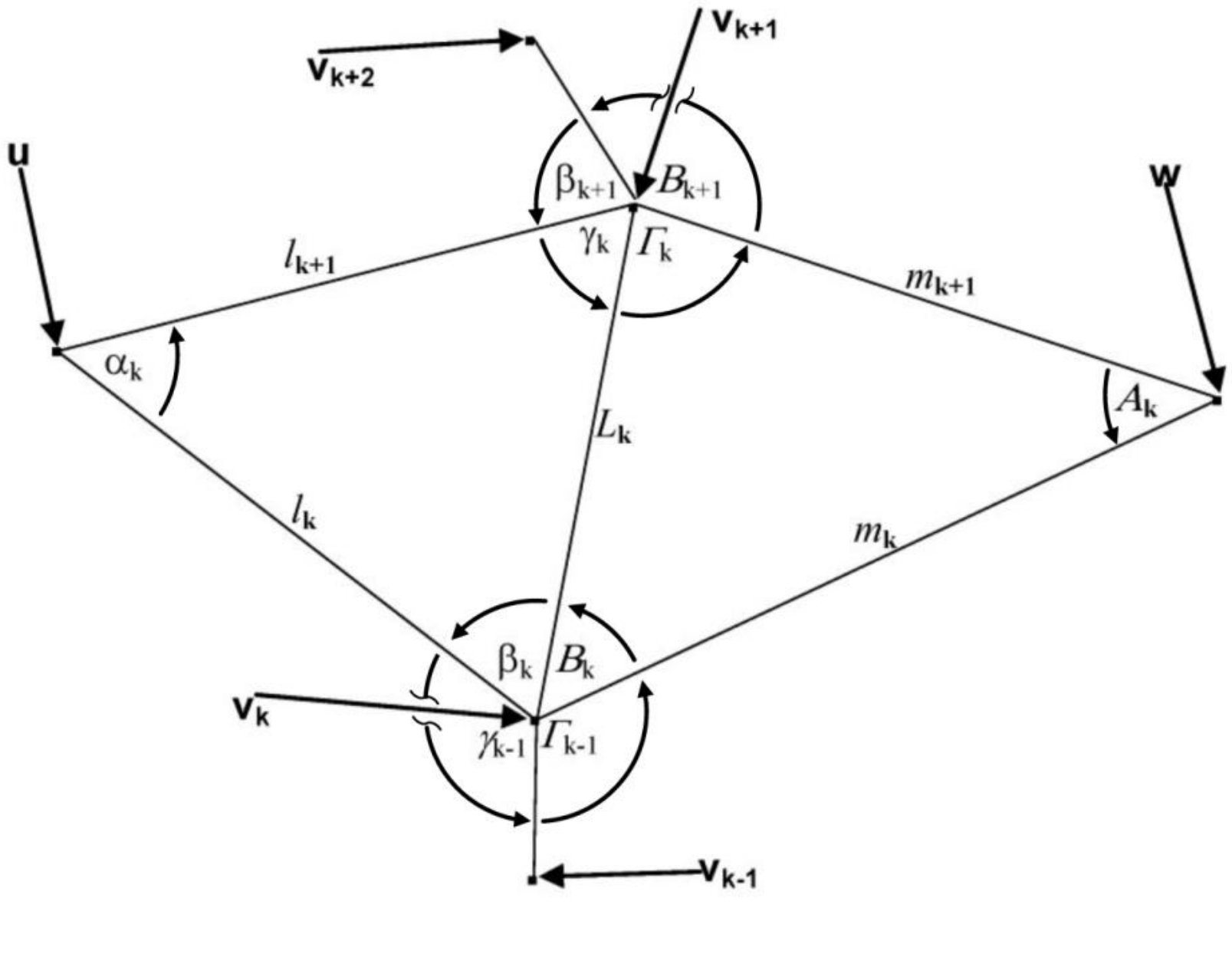

**Figure 2.** Notation.

Fig. 2 illustrates two typical faces, $f_k$ and $F_k$, $\mathbf{uv_{k+1}v_k}$ and $\mathbf{wv_kv_{k+1}}$ respectively, and is by convention the view of the outsides of the faces. The vectors $\mathbf{(v_{k+1}\text{-}u)X(v_k\text{-}u)}$ and $\mathbf{(v_k\text{-}w)X(v_{k+1}\text{-}w)}$, where $\mathbf{X}$ is the vector cross product operator, determine the direction vectors associated with these faces. They point into the page.

## 4. I-OEE Generalization – the first Bricard type

Bricard octahedra of the first type [1, Eq. 11] have six pairs of edges that are of equal length and are arranged so that opposite edges have equal lengths. There are four pairs of identical faces. The two apical vertexes $\mathbf{u}$ and $\mathbf{w}$ have four identical faces that are ordered identically but offset by two faces at the base of the caps. The faces that are identical on the caps do not share a common edge but are opposite of one another. Edge lengths associated with the vertexes $\mathbf{u}$ and $\mathbf{w}$ have the following relationships: $m_1=l_3$, $m_2=l_4$, $m_3=l_1$, $m_4=l_2$, while at the base of the cap the edge lengths are OEE: $L_3=L_1$ and $L_4=L_2$.

These relationships can be achieved in general when N=2M>4 by treating 3M edge lengths as independent parameters and restricting the edge lengths as follows:

$$m_k=l_{k+M},\quad m_{k+M}=l_k \text{ and } L_{k+M}=L_k \text{ for } k=1..M. \tag{1}$$

These relationships are shown schematically in Fig. 3 where the lengths of 3M independent parameters are displayed in their relative positions on the suspension. The top line represents the edge lengths on the $\mathbf{u}$ cap, the middle line the base of the cap (edge lengths and vertexes) and the bottom line the edge lengths on the $\mathbf{w}$ cap.

| | | | | | | | | | | | | | | | | | | |
|---|---|---|---|---|---|---|---|---|---|---|---|---|---|---|---|---|---|---|
| $l_1$ | | $l_2$ | | $l_3$ | | $\ldots l_{M-1}$ | | $l_M$ | | $l_{M+1}$ | | $l_{M+2}$ | | $\ldots l_{N-1}$ | | $l_N$ | | $l_1$ |
| $\mathbf{v_1}$ | $L_1$ | $\mathbf{v_2}$ | $L_2$ | $\mathbf{v_3}$ | $L_3$ | $\ldots\mathbf{v_{M-1}}$ | $L_{M-1}$ | $\mathbf{v_M}$ | $L_M$ | $\mathbf{v_{M+1}}$ | $L_1$ | $\mathbf{v_{M+2}}$ | $L_2$ | $\ldots\mathbf{v_{N-1}}$ | $L_{M-1}$ | $\mathbf{v_N}$ | $L_M$ | $\mathbf{v_1}$ |
| $l_{M+1}$ | | $l_{M+2}$ | | $l_{M+3}$ | | $\ldots l_{N-1}$ | | $l_N$ | | $l_1$ | | $l_2$ | | $\ldots l_{M-1}$ | | $l_M$ | | $l_{M+1}$ |

**Figure 3.** I-OEE Edge Length Relationships.

Additionally we associate the I-OEE sub-type with the coordinate model:

$\mathbf{u}(0,0,z)$, $\mathbf{w}(0,0,-z)$, $\{\mathbf{v_k}(x_k,y_k,z_k)$, $\mathbf{v_{k+M}}(-x_k,y_k,-z_k)$, k=1..M$\}$.

By direct evaluation of the edge lengths it is easy to verify that this coordinate model and the edge length relationships (Eqs. 1) are consistent. Assuming flexibility it is seen that the characteristic axial motion of Bricard octahedra of the first type is preserved in this generalization. The two apical vertexes, $\mathbf{u}$ and $\mathbf{w}$, will move in a line while at the base of the caps, pairs of non-apical vertexes will move symmetrically about an axis that is orthogonal to this line.

From this coordinate model it can be seen that the dihedral angles of edges associated with symmetrically opposite vertexes $\mathbf{v_k}$ and $\mathbf{v_{k+M}}$ for k=1.M and vertexes $\mathbf{u}$ and $\mathbf{w}$ are conjugate in value. Specifically the following relationships hold for dihedral angles:

$$\delta_k = 2\pi - \Delta_{k+M},\quad \delta_{k+M} = 2\pi - \Delta_k \text{ and } \varepsilon_k = 2\pi - \varepsilon_{k+M} \text{ for } k=1..M. \tag{2}$$

## 5. II-AEE Generalization – the second Bricard type

Bricard octahedron of the second type [1, Eq. 18] have six pairs of edges that are of equal length and four pairs of identical faces that are arranged differently than in the first Bricard type. In the II-AEE characterization each of the two apical vertexes have four distinct faces that are arranged so that on the base of the caps two adjacent edges are of equal length. Edge lengths associated with the caps at **u** and **w** have the following relationships: $m_1=l_1$, $m_2=l_4$, $m_3=l_3$, $m_4=l_2$. At the base of the cap there are edge lengths that are adjacent and of equal length, ie. are AEE: $L_3=L_2$, $L_4=L_1$. (The choice of which edges on the cap base are assigned to have equal lengths is by convention.)

These relationships can be generalized by treating 3M edge lengths as independent parameters and restricting the edge lengths as follows:

$$m_1=l_1 \text{ and}$$
$$m_k=l_{N-k+2} \text{ for k=2..N and} \qquad (3)$$
$$L_{N+1-k}=L_k \text{ for k=1..M.}$$

An arrangement of the 3M independent edge length parameters is illustrated in Fig. 4.

$l_1$ $l_2$ $l_3$ $\ldots l_{M-1}$ $l_M$ $l_{M+1}$ $l_{M+2}$ $\ldots l_{N-1}$ $l_N$ $l_1$
$\mathbf{v_1}$ $L_1$ $\mathbf{v_2}$ $L_2$ $\mathbf{v_3}$ $L_3 \ldots \mathbf{v_{M-1}}$ $L_{M-1}$ $\mathbf{v_M}$ $L_M$ $\mathbf{v_{M+1}}$ $L_M$ $\mathbf{v_{M+2}}$ $L_{M-1}$ $\ldots \mathbf{v_{N-1}}$ $L_2$ $\mathbf{v_N}$ $L_1$ $\mathbf{v_1}$
$l_1$ $l_N$ $l_{N-1}$ $\ldots l_{M+3}$ $l_{M+2}$ $l_{M+1}$ $l_M$ $\ldots l_3$ $l_2$ $l_1$

**Figure 4.** II-AEE Edge Length Relationships.

A coordinate model for the II-AEE sub-type that is consistent with the edge length relationships (Eqs. 3) is the following:

$$\mathbf{u}(0,0,z),\ \mathbf{w}(0,0,-z),\ \mathbf{v_1}(x_1,y_1,0),\ \mathbf{v_{M+1}}(x_{M+1},y_{M+1},0),$$
$$\{\mathbf{v_k}(x_k,y_k,z_k),\ \mathbf{v_{M+k}}(x_{M-k+2},y_{M-k+2},-z_{M-k+2}),\ k=2..M\}.$$

Assuming flexibility, the characteristic planar motion of the second type of Bricard octahedra is exhibited in this generalization. The two vertexes, $\mathbf{v_1}$ and $\mathbf{v_{M+1}}$, will move in a plane while the remaining vertexes move symmetrically with respect to this plane.

From this coordinate model it can be seen that the dihedral angles of edges associated with asymmetrically opposite vertexes $\mathbf{v_k}$ and $\mathbf{v_{N-k+1}}$ for k=1.M and vertexes **u** and **w** are conjugate in value. Specifically the following relations hold for dihedral angles:

$$\delta_1 = 2\pi - \Delta_1,$$
$$\delta_k = 2\pi - \Delta_{N-k+2}, \text{ for k=2..N,} \qquad (4)$$
$$\varepsilon_k = 2\pi - \varepsilon_{N-k+1} \text{ for k=1..M.}$$

## 6. II-OEE Generalization – the second Bricard type

In the II-OEE characterization of Bricard octahedron of the second type, the two caps at vertexes **u** and **w** have different faces and each cap has pairs of identical faces that are arranged opposite of one another on the cap; consequently the I-OEE arrangement holds for the edges on the base of the cap with an alternating arrangement of edge lengths at each apical vertex. Edge lengths associated with the caps at **u** and **w** have the following relationships: $l_3=l_1$, $l_4=l_2$, $m_3=m_1$, $m_4=m_2$, $L_3=L_1$, $L_4=L_2$. In general this leads to M independent edge lengths on each cap and M common edge lengths at the base:

$$l_{k+M} = l_k,$$
$$m_{k+M} = m_k \text{ and}$$
$$L_{k+M} = L_k \text{ for } k=1..M. \qquad (5)$$

These relationships are illustrated in Fig. 5.

$$l_1 \quad l_2 \quad l_3 \quad \ldots l_{M-1} \quad l_M \quad l_1 \quad l_2 \quad \ldots l_{M-1} \quad l_M \quad l_1$$
$$\mathbf{v_1}\ L_1\ \mathbf{v_2}\ L_2\ \mathbf{v_3}\ L_3 \ldots \mathbf{v_{M-1}}\ L_{M-1}\ \mathbf{v_M}\ L_M\ \mathbf{v_{M+1}}\ L_1\ \mathbf{v_{M+2}}\ L_2 \ldots \mathbf{v_{N-1}}\ L_{M-1}\ \mathbf{v_N}\ L_M\ \mathbf{v_1}$$
$$m_1 \quad m_2 \quad m_3 \quad \ldots m_{M-1} \quad m_M \quad m_1 \quad m_2 \quad \ldots m_{M-1} \quad m_M \quad m_1$$

**Figure 5.** II-OEE Edge Length Relationships.

A coordinate model for the II-OEE sub-type that is consistent with these edge length relationships (Eqs. 5) is the following:

$$\mathbf{u}(x,y,0),\ \mathbf{w}(X,Y,0),\ \{\mathbf{v_k}(x_k,y_k,z_k),\ \mathbf{v_{k+M}}(x_k,y_k,-z_k),\ k=1..M\} \text{ with } x_1=y_1=0.$$

Assuming flexibility, this generalization also exhibits a form of planar symmetry. Non-apical vertexes will move symmetrically with respect to the plane of motion of the apical vertexes.

From this coordinate model it can be seen that the dihedral angles of edges associated with symmetrically opposite vertexes $\mathbf{v_k}$ and $\mathbf{v_{k+M}}$ for k=1.M and vertexes $\mathbf{u}$ and $\mathbf{w}$ are seen to be conjugate in value. Specifically the following relations hold for dihedral angles:

$$\delta_k = 2\pi - \delta_{k+M},$$
$$\Delta_k = 2\pi - \Delta_{k+M}, \qquad (6)$$
$$\varepsilon_k = 2\pi - \varepsilon_{k+M} \text{ for } k=1..M.$$

## 7. The third Bricard type

A primary characteristic of Bricard octahedra of the third type [1, Sec. 11] that we seek to retain in our generalization is the ability of these octahedra to have two positions in which all vertexes are co-planar or alternately stated, to have two positions in which all dihedral angles simultaneously have the values $0^o$ or $180^o$. In Bricard octahedra this type of folding is enabled by vertexes that either have opposite angles that are equal or that are supplementary in value (OAE or OAE) and specifically by having four OAE vertexes and two OAS vertexes that are opposite of one another [7, Fig. 5].

For Bricard octahedra of the third type the III-OAE designation implies that the apical vertexes $\mathbf{u}$ and $\mathbf{w}$ are OAE while on the cap base there are two OAS vertexes, $\mathbf{v_1}$ and $\mathbf{v_3}$ and two OAE vertexes, $\mathbf{v_2}$ and $\mathbf{v_4}$. (The choice of which pair of vertexes is treated as OAS is by convention.) With the III-OAS designation the two apical vertexes, $\mathbf{u}$ and $\mathbf{w}$, are OAS and there are four OAE vertexes on the base of the cap. Further the folding that is possible with these octahedra can be characterized by two simple equations that define the relationships between the face angles at the caps. Equations for the angles at the apical vertex $\mathbf{w}$ are not shown since they are identical in form.

For the III-OAE sub-type:

$$\alpha_1 + \alpha_2 - \alpha_3 - \alpha_4 = 0 \text{ and}$$
$$\alpha_1 - \alpha_2 - \alpha_3 + \alpha_4 = 0.$$

These define an open co-planar folding in which there are two sets of equal and opposite dihedral angles. For the first: $\delta_1 = \delta_3 = 0$ and $\delta_2 = \delta_4 = \pi$, while for the second: $\delta_2 = \delta_4 = 0$ and $\delta_1 = \delta_3 = \pi$. These equations simplify to

$\alpha_1 = \alpha_3$ and
$\alpha_2 = \alpha_4$.

For the III-OAS sub-type:

$\alpha_1 + \alpha_2 + \alpha_3 + \alpha_4 = 2\pi$ and
$\alpha_1 - \alpha_2 + \alpha_3 - \alpha_4 = 0$.

The first defines an open co-planar folding in which $\delta_1=\delta_2=\delta_3=\delta_4=\pi$. The second defines a compact co-planar folding in which $\delta_1=\delta_2=\delta_3=\delta_4=0$. These equations simplify to

$\alpha_1 = \pi-\alpha_3$ and
$\alpha_2 = \pi-\alpha_4$.

The above equations do not characterize the octahedra as completely as with sub-types I-OEE, II-AEE and II-OEE polyhedra where edge lengths are given explicitly and face angles readily angles readily determined from law of sines relationships. To achieve this level of characterization we parameterize two adjacent faces of the octahedron with the parameter set $\{l_1,\alpha_1,\beta_1,\alpha_2,\beta_2\}$ (Bricard octahedra of the third type are known [8] to have five parameters.) and use law of sines relationships and relationships that exist between the dihedral angles of opposite vertexes to compute edge lengths and face angles. Computational details are provided in Appendix A.

In addition to Bricard octahedra of the third type, the flexible dodecahedral suspension that is described in [5] provides further information about how a generalization of this type of suspension may proceed. Inspection of this polyhedron reveals that there are two OAS vertexes, $\mathbf{v_1}$ and $\mathbf{v_3}$ say, and four OAE vertexes on the cap base and that the configuration of cap vertex face angles that describe the flat co-planar folding is given (using the notation of this paper) by:

$\alpha_1 + \alpha_2 - \alpha_3 - \alpha_4 - \alpha_5 - \alpha_6 = 0$ and
$\alpha_1 - \alpha_2 - \alpha_3 + \alpha_4 - \alpha_5 + \alpha_6 = 0$.

These equations simplify to

$\alpha_1 = \alpha_3 + \alpha_5$ and
$\alpha_2 = \alpha_4 + \alpha_6$.

## 8. III-OAE Generalization – the third Bricard type

The III-OAE generalization that is suggested by the above discussion is one in which there are two OAS vertexes, $\mathbf{v_1}$ and $\mathbf{v_L}$ for some L where $3 \leq L \leq (N-2)$ along with N-2 OAE vertexes on the cap base and a "generalized" OAE cap with the following configuration of cap vertex face angles equations:

$$\sum_{k=1..L-1} \alpha_k - \sum_{k=L..N} \alpha_k = 0 \text{ and}$$

$$\sum_{k=1..L-1} (-)^{k+1}\alpha_k - \sum_{k=L..N} (-)^{k+1}\alpha_k = 0. \qquad (7)$$

The first of these equations defines an open co-planar folding with $\delta_1=\delta_L=0$ and $\delta_k=\pi$ for k=2..L-1 and k=L+1..N. The second represents a co-planar folding in which $\delta_1=\delta_L=\pi$ and $\delta_k=0$ for k=2..L-1 and k=L+1..N.

Equations (7) simplify to:

$$\sum_{k=1..k_x} \alpha_{2k-1} = \sum_{k=k_x+1..M} \alpha_{2k-1} \text{ and}$$

$$\sum_{k=1..k_x-1} \alpha_{2k} = \sum_{k=k_x..M} \alpha_{2k} \quad \text{when L is even } (L = 2k_x) \text{ or}$$

$$\sum_{k=1..k_x} \alpha_{2k} = \sum_{k=k_x+1..M} \alpha_{2k} \quad \text{when L is odd } (L = 2k_x+1).$$

With this configuration of face angles at vertex **u** it is easy to show that the face angles at vertex **w**, angles $A_k$ for k=1..N also satisfy equations (7). These angles can be computed directly based upon the OAE and OAS assignments to the vertexes to yield:

$$\begin{aligned}
A_1 &= \pi - \alpha_N - \beta_N - \beta_2, \\
A_2 &= \alpha_1 + \beta_1 - \beta_3, \\
A_3 &= \alpha_2 + \beta_2 - \beta_4, \\
&\;\;. \\
&\;\;. \\
A_k &= \alpha_{k-1} + \beta_{k-1} - \beta_{k+1}, \text{ k=2..L-2 and k=L+1..N-1}, \\
&\;\;. \\
&\;\;. \\
A_{L-1} &= -\pi + \alpha_{L-2} + \beta_{L-2} + \beta_L, \\
A_L &= \pi - \alpha_{L-1} - \beta_{L-1} - \beta_{L+1}, \\
&\;\;. \\
&\;\;. \\
A_N &= -\pi + \alpha_{N-1} + \beta_{N-1} + \beta_1.
\end{aligned}$$

Substitution of these values into equations (7) yields:

$$\sum_{k=1..L-1} A_k - \sum_{k=L..N} A_k = \sum_{k=1..L-1} \alpha_k - \sum_{k=L..N} \alpha_k = 0 \text{ and}$$

$$\sum_{k=1..L-1} (-)^{k+1} A_k - \sum_{k=L..N} (-)^{k+1} A_k = . \sum_{k=1..L-1} (-)^{k+1} \alpha_k - \sum_{k=L..N} (-)^{k+1} \alpha_k = 0$$

As with Bricard octahedra of the third type the above characterization (Eqs. 7) is not sufficient to fully parameterize the suspension of interest. To accomplish a complete parameterization of all faces we use the parameter set $\{l_1,\alpha_1,\beta_1,\alpha_2,\beta_2,\alpha_3,\alpha_5,\ldots,\alpha_{N-3}\}$ and employ a recursive strategy that is based upon relationships between dihedral angles associated with the vertex pairs:

$(\mathbf{v_1},\mathbf{v_3}),(\mathbf{v_2},\mathbf{v_5}), (\mathbf{v_4},\mathbf{v_7}),\ldots. (\mathbf{v_{2k}},\mathbf{v_{2k+3}}),\ldots (\mathbf{v_{N-4}},\mathbf{v_{N-1}}), (\mathbf{v_{N-2}},\mathbf{v_N})$ for k=1..M-2.

This is a generalization of vertex pairs that are observed in Bricard octahedra of the third type and in the flexible dodecahedral suspension [5] discussed earlier. Appendix B contains the details of the parameterization computations.

This generalization also leads to equal or conjugate relationships for the dihedral angles associated with the vertex pairs shown above. Specifically the following relations hold for dihedral angles:

$$\delta_k = \Delta_j \mid \pi - \Delta_j \mid 2\pi - \Delta_j \mid 3\pi - \Delta_j,$$
$$\Delta_k = \delta_j \mid \pi - \delta_j \mid 2\pi - \delta_j \mid 3\pi - \delta_j \mid, \quad (8)$$
$$\varepsilon_k = \varepsilon_j \mid 2\pi - \varepsilon_j \text{ for all (k, j) pairs as shown above.}$$

## 9. III-OAS Generalization – the third Bricard type

With the III-OAS designation, flexible octahedra have apical vertexes **u** and **w** that are OAS while all non-apical vertexes are OAE. Generalization is made in a straightforward manner by designating the N non-apical vertexes as OAE, restricting one of the cap face angle configurations to angles that define an open folding in which all dihedral angles $\delta_k=\pi$ and restricting a second to a folding that is general in nature and represents a compact folding in which all faces are folded and stacked together with all dihedrals $\delta_k=0$ for k=1..N:

$$\sum_{k=1..N} \alpha_k = 2K\pi \text{ and}$$
$$\sum_{k=1..N} (-)^{k+1}\alpha_k = 0 \quad (9)$$

for some integer K>0. These equations simplify to:

$$\sum_{k=1..M} \alpha_{2k-1} = K\pi \text{ and}$$

$$\sum_{k=1..M} \alpha_{2k} = K\pi.$$

With this characterization it is easy to see that the above equations are also satisfied by the face angles at vertex **w**. To complete the characterization of this sub-type we employ the same parameter set and approach as used for sub-type III-OAE (Appendix B).

## 10. Construction of Suspensions

The generalizations that are defined in Sections 4 – 9 provide or lead to a complete description of polyhedral suspensions in terms of edge lengths and face angles. However, they do not contain explicit information as to how the polyhedra can be constructed or flexion achieved or even if construction is possible. Here we describe recursive approaches that address the construction problem and the problem of testing the flexibility or rigidity of suspensions constructed with specific parameterizations. While each step of these recursions is relatively simple they are regarded as forming the basis for computer implementations.

For the purpose of construction the dihedral angle $\varepsilon_1$ on the edge $\mathbf{v_1v_2}$ is treated as the variable of flexion. To initiate the recursion the two faces $\mathbf{uv_2v_1}$ and $\mathbf{wv_1v_2}$ are positioned by fixing the face $\mathbf{uv_2v_1}$ at some convenient coordinate location followed by rotation about $\mathbf{v_2v_1}$ by $\varepsilon_1$ to locate **w**, eg. $\mathbf{u}(l_1\sin\beta_1,0,-l_1\cos\beta_1)$, $\mathbf{v_1}(0,0,0)$, $\mathbf{v_2}(0,0,-L_1)$, and $\mathbf{w}(m_1\sin B_1\cos\varepsilon_1, m_1\sin B_1\sin\varepsilon_1, -m_1\cos B_1)$.

Subsequent recursive construction is based upon the ability to compute the dihedral angles at a tetrahedral vertex given that one dihedral angle and all face angles at the vertex are known. We proceed as follows: at the k-th step, for k=2..N-1, the three dihedral angles $\{\delta_k, \Delta_k, \varepsilon_k\}$ at vertex $\mathbf{v_k}$ are determined from a solution of equation (C-5) as applied to vertex $\mathbf{v_k}$ and using the value of $\varepsilon_{k-1}$ that was determined at vertex $\mathbf{v_{k-1}}$. From the dihedral angles determined at vertex $\mathbf{v_k}$ the position of vertex $\mathbf{v_{k+1}}$ is determined by rotation about edge $\mathbf{uv_k}$ by $\delta_k$. In this manner the two faces $\mathbf{uv_{k+1}v_k}$ and $\mathbf{wv_kv_{k+1}}$ are defined at each stage of the recursion. The two faces $\mathbf{uv_1v_N}$ and $\mathbf{wv_Nv_1}$ are defined at termination of the recursion based upon the position of vertex $\mathbf{v_N}$.

This approach is over-determined in the sense that some of the edge lengths of the suspension are not explicitly used to position vertexes. Specifically the edge length $l_N=|\mathbf{v_Nv_1}|$ is not used. Since this edge length is indeterminate the construction sequence may or may not result in the construction of a polyhedron with the specified edge length and face angle parameters. This situation is exacerbated by the fact that the solution of equation (C-5) is quadratic in nature; consequently there are $2^{N-2}$ potential constructions of interest.

Fortunately for sub-types I-OEE, II-AEE and II-OEE the symmetric relations defined by the coordinate models associated with Figs. 3, 4 and 5 can be used to complete a specific construction given that the positions of vertexes $\mathbf{v_1}$ through $\mathbf{v_{M+1}}$ have been determined by the above approach. For all vertexes with an index >M+1 the relevant solution of equation (C-5) is simply the other solution associated with its symmetric vertex. This leads to $2^{M-1}$ unique constructions that are flexible. For the sub-types III-OAE and III-OAS the ambiguity in the choice of dihedral angle at particular vertexes can be resolved using the sequence of partial constructions that was used in the computation of the complete parameter set as discussed in Appendix B.

Once construction is complete vertex positions are known and all dihedral angles and their derivatives can be determined. From these a straightforward test for flexibility that is both necessary and sufficient can be implemented. Vectors $\mathbf{u'}(0,0,0)$, $\mathbf{w'}(-m_1\sin\varepsilon_1\sin B_1, 0, m_1\cos\varepsilon_1\sin B_1)$, $\mathbf{v'}_1(0,0,0)$ and $\mathbf{v'}_2(0,0,0)$ where $(\ldots)'$ denotes differentiation with respect to $\varepsilon_1$ are treated as initial conditions. Recursively, for k=3..N, $\mathbf{v_k'}$ and the derivative of dihedral angles at $\mathbf{v_k}$ are determined from the derivative of the positions of $\mathbf{u'}$ and $\mathbf{v'_{k-1}}$ and the dihedral angle derivatives at $\mathbf{v_{k-1}}$. Dihedral angle derivatives are computed from equations (C-4), (C-5) and (C-6). This enables the evaluation of the derivative term:

$$\frac{d|\mathbf{p}\text{-}\mathbf{q}|^2}{d\varepsilon_1} = 2|\mathbf{p}\text{-}\mathbf{q}| \cdot \frac{d|\mathbf{p}\text{-}\mathbf{q}|}{d\varepsilon_1} \tag{10}$$

for any edge $\mathbf{pq}$ of the polyhedron of interest and in particular for $\mathbf{p}=\mathbf{v_N}$ and $\mathbf{q}=\mathbf{v_1}$. The vanishing of this term serves as the basis for testing the validity of the construction.

## 11. Flexible Suspension Constructions

In this section we discuss example flexible suspensions that have been constructed using the generalizations that are defined in Sections 4 through 9 and the construction methods defined in Section 10. These examples have been constructed as a validation of the generalizations and methods of construction and to provide a better understanding of the details of flexible suspensions. Examples are provided for all five sub-types and for N=4,6,8,10,12,14,16. Example octahedra have been included for the sake of completeness.

A detailed description of the examples, the parameter values used for construction and other details of interest, is found in Appendix D. Sufficient data is provided there to support the construction of any of the examples. Here we provide some general observations and comments. All parameter values used for these example constructions should be regarded as nominal values that have been determined to illustrate the suspension of interest and generally do not have any particular significance although in practice some experimentation may have been employed to ensure that each parameter value is a whole number. In all cases the value of parameters can be varied, individually or in concert, over some interval to produce similar example flexible polyhedra.

Sub-type I-OEE, II-AEE and II-OEE examples are created using the parameter set $\{l_{k,}$ k=1..N; $L_k$, k=1..M$\}$, constructed as defined in Section 10 and then transformed into the coordinate models associated with Figs. 3, 4 and 5 respectively for the purpose of display. The parameterization of suspensions of these sub-types is a straightforward exercise. Any collection of triangular faces that are constructed to conform to the specifications of equations (1, 3 or 5) can be assembled into $2^{M-1}$ flexible polyhedral suspensions. III-OAE and III-OAS examples are created with the parameter set $\{l_1,\alpha_1,\beta_1,\alpha_2,\beta_2,\alpha_3,\alpha_5,\ldots\alpha_{2k+1},\ldots\alpha_{N-3}$; k=1..M-2$\}$and constructed as defined in Section 10. For these sub-types, the determination of valid parameters is the primary computational challenge associated with these suspensions and has been accomplished by manual experimentation with automated evaluation of the multiplicity of choices as described in Appendix B.

All examples have been validated by a number of computations that are evaluated over the full range of flexion. Some computations are general in nature and some of which are sub-type specific. Those of a general nature are the vanishing of the derivative term defined by equation (10), zero value of the oriented volume and constant value of the total mean curvature. Sub-type specific terms are the relationships between the various dihedral angles as defined by equations ( 2, 4, 6 and 8). Additionally there exist a number of relationships between the solid angle of vertex pairs that are constant during flexion. Typically the solid angles of each of the vertex pairs of interest add together to $=4\pi$.

## 12. Conclusion

In this paper we have described five families of polyhedral suspensions that can be inferred from generalizations of the description of the geometric characteristics of Bricard octahedra. In the case of the sub-type III-OAE a new, recently reported [5] flexible dodecahedral suspension is also used as a basis for generalization. All generalizations reduce to Bricard octahedra characterizations for N=4.

A description of the parameterization and complete characterization of the suspensions along with construction methods and a numerical test for flexibility that is necessary and sufficient are included. From these we have created constructions of flexible polyhedra that are representative of all five families. Consequently, we form the following conjectures that flexible suspensions of indefinite size exist that conform to the generalizations.

**Conjecture I**. For any I-OEE, II-AEE or II-OEE suspension, as defined in Secs. 4, 5 and 6 respectively, there exists a finite interval $[\varepsilon_m, \varepsilon_x]$ over which the variable $\varepsilon_1$ can be varied and construction of the suspension achieved.

In other words the descriptions given for these suspensions are necessary and sufficient conditions for flexibility.

The case of the sub-types III-OAE and III-OAS is complicated by the fact that the parameter sets are insufficient to construct the suspension and the recursive computations described in Appendix B are required to generate the necessary edge lengths and face angles. Also not all parameter sets lead to realizable values. However, in a number of examples all parameter sets that led to realizable values also yielded flexible suspensions. Consequently we form the following:

**Conjecture II**. Any III-OAE or III-OAS parameter set that leads to realizable edge lengths and face angles under the computations as described in Appendix B define a flexible suspension.

We close by noting that in addition to proving the above conjectures for all even N>4 there exist a number of additional related studies that may be of interest; chief among them is the description of flexible suspensions for odd N≥5. Also it is relatively straight forward exercise to extend the flexible suspensions described in this paper into larger flexible polyhedra of both genus 0 and 1 and of indefinite size by generalizing an approach [6] taken to extend Bricard octahedra. Finally, it is noted that the descriptions are not exhaustive. In the course of preparing this paper we were able to construct a flexible suspension with N=8 in which non-apical vertexes were alternately OAE and OAS but we were unable to form a generalization for N>8. Also in this regard there are numerous other foldings analogous to equation (7) that could potentially yield flexible suspensions that exhibit two flat positions.

## Acknowledgements

The author wishes to thank his wife Verla Nelson for her encouragement and support.

**Gerald Nelson** is a retired software engineer; was employed by Honeywell, Inc. (when it was a Minnesota based company) and by MTS Systems Corp. of Eden Prairie, Minnesota. He has a Masters degree in Mathematics from the University of Minnesota. His interest in flexible polyhedra stems from reading an article [11] about the Bellows Conjecture.

References.

## Appendix A.

This appendix describes an approach that is used to fully parameterize Bricard octahedra of the third type from the parameter set $\{l_1,\alpha_1,\beta_1,\alpha_2,\beta_2\}$ and as will be seen in Appendix B, is used as an initial step in fully parameterizing suspensions for N>4. This parameter set completely defines the two adjacent faces $\mathbf{uv_2v_1}$ and $\mathbf{uv_3v_2}$ but parameters on the other faces are unresolved. The situation is shown schematically in Fig. A-1 which depicts parameters from the parameter set, e.g. $\boldsymbol{\alpha_1}$, parameters that are readily computed from the OAS and OAE assignments and the law of sines, e.g. $\alpha_3$, and the parameters that are unresolved and determined as described in this Appendix, e.g. $A_1$. It is apparent from inspection that the value of one parameter, say $\beta_3$, is sufficient to complete the parameterization.

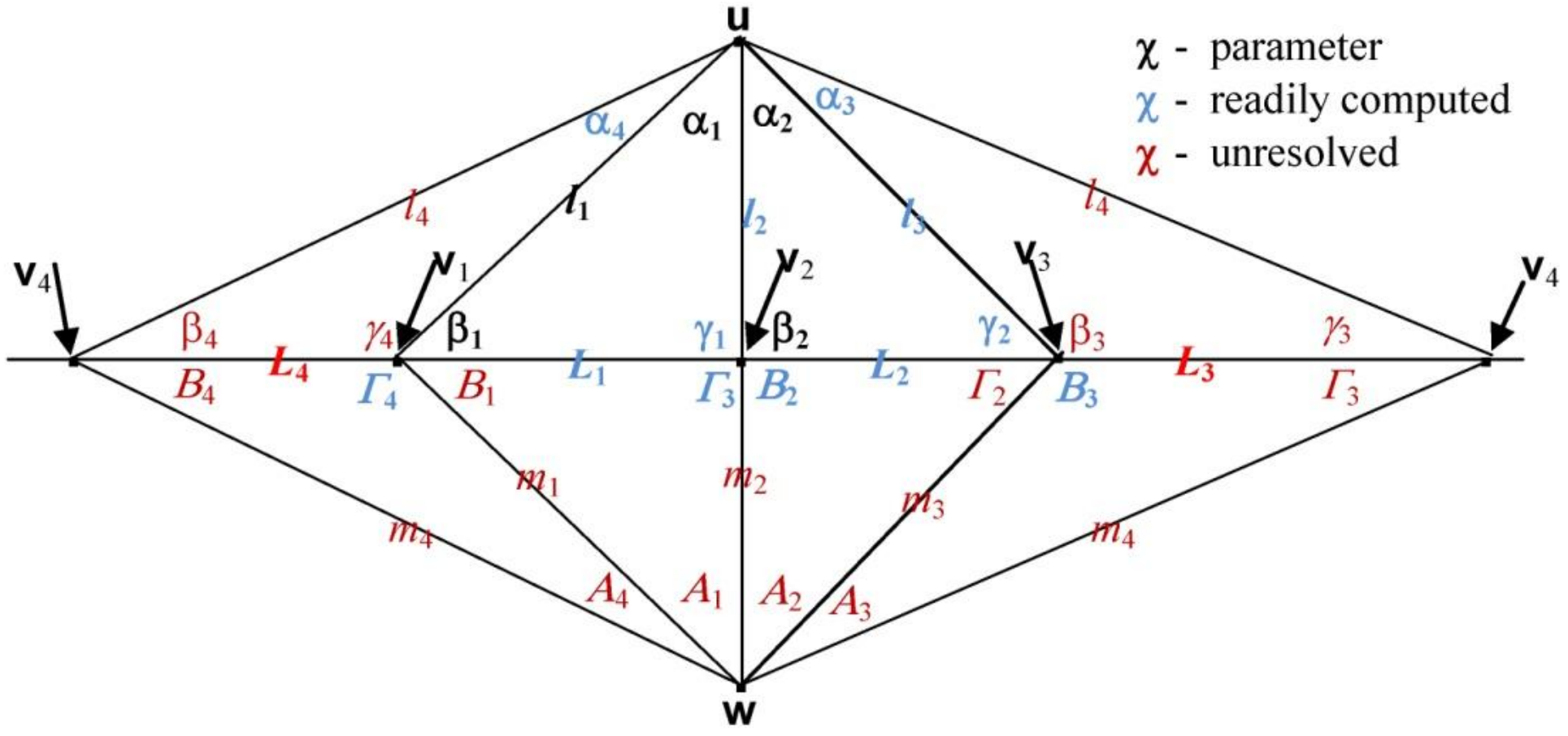

**Figure A-1.** Octahedron Parameterization.

Determination of $\beta_3$ proceeds as follows. By applying the law of sines to the faces $\mathbf{wv_1v_2}$ and $\mathbf{wv_2v_3}$ and eliminating the edge length $m_2$ from the resulting equations it is easy to see that the two angles $\Gamma_2$ and $B_1$ are related by:

$$\cot\Gamma_2 = a \cot B_1 + b, \qquad \text{(A-1)}$$

where

$$a = \frac{L_2 \sin\Gamma_1}{L_1 \sin B_2} \quad \text{and}$$

$$b = \frac{L_2 \cos\Gamma_1 - L_1 \cos B_2}{L_1 \sin B_2}.$$

When $\mathbf{v_3}$ is OAE $\Gamma_2=\beta_3$ and equation (A-1) can be written as

$$\cot\beta_3 = a \cot B_1 + b. \qquad \text{(A-2)}$$

When $\mathbf{v_3}$ is OAS $\Gamma_2=\pi-\beta_3$ and equation (A-1) becomes

$$\cot\beta_3 = -a \cot B_1 - b. \qquad \text{(A-3)}$$

A second relationship that exists between the variables $\beta_3$ and $B_1$ can be derived from a consideration of the equations [1, Eqs. 4 and 5] that describe the deformation of tetrahedral vertexes in Bricard octahedra of the third type. When these equations are applied to the vertex $\mathbf{v_j}$, one of the following equations that relate the dihedral angles $\varepsilon_j$ and $\Delta_j$ and the face angles $\beta_j$ and $B_j$ must be true:

$\mathbf{v_j}$ is OAE: $t_h(\varepsilon_j)t_h(\Delta_j) = b_1(\beta_j, \mathcal{B}_j)$ or $= b_2(\beta_j, \mathcal{B}_j)$ and (A-4)

$\mathbf{v_j}$ is OAS: $t_h(\varepsilon_j)c_h(\Delta_j) = -b_1(\beta_j, \mathcal{B}_j)$ or $= -b_2(\beta_j, \mathcal{B}_j)$,

where the functions $b_1$, $b_2$, $t_h$ and $c_h$ are defined by:

$$b_1(\beta, B) = \frac{\sin\left(\frac{\beta - B}{2}\right)}{\sin\left(\frac{\beta + B}{2}\right)},$$

$$b_2(\beta, B) = \frac{\cos\left(\frac{\beta - B}{2}\right)}{\cos\left(\frac{\beta + B}{2}\right)}, \qquad \text{(A-5)}$$

$$t_h(\delta) = \tan\left(\frac{\delta}{2}\right) \text{ and}$$

$$c_h(\delta) = \cot\left(\frac{\delta}{2}\right).$$

In Bricard octahedra of the third type several relationships between dihedral angles exist as a consequence of equation (C-6):

$\cos\varepsilon_1 = \cos\varepsilon_2 = \cos\varepsilon_3 = \cos\varepsilon_4$,

$\cos\delta_1 = \cos\delta_3 = \cos\Delta_1 = \cos\Delta_3$ and

$\cos\delta_2 = \cos\delta_4 = \cos\Delta_2 = \cos\Delta_4$.

Thus each of the dihedral angles in these equations are equal or conjugate in value and from equation (A-4) one of the following is true:

$b_1(\beta_j, \mathcal{B}_j) = \pm b_1(\beta_k, \mathcal{B}_k)$, (a,b)

$b_1(\beta_j, \mathcal{B}_j) = \pm b_2(\beta_k, \mathcal{B}_k)$, (c,d) (A-6)

$b_2(\beta_j, \mathcal{B}_j) = \pm b_1(\beta_k, \mathcal{B}_k)$ or (e,f)

$b_2(\beta_j, \mathcal{B}_j) = \pm b_2(\beta_k, \mathcal{B}_k)$. (g,h)

for the vertex pair ($\mathbf{v_j}$,$\mathbf{v_k}$) when j=1 and k=3. Using well known trigonometric relationships equations (A-6) reduce to:

$1 = t_h(\mathcal{B}_j)\, c_h(\beta_j)\, c_h(\mathcal{B}_k)\, t_h(\beta_k)$, (a)

$1 = t_h(\mathcal{B}_j)\, c_h(\beta_j)\, t_h(\mathcal{B}_k)\, c_h(\beta_k)$, (b)

$1 = -t_h(\mathcal{B}_j)\, c_h(\beta_j)\, t_h(\mathcal{B}_k)\, t_h(\beta_k)$, (c)

$1 = -t_h(\mathcal{B}_j)\, c_h(\beta_j)\, c_h(\mathcal{B}_k)\, c_h(\beta_k)$, (d) (A-7)

$1 = -c_h(B_j)\ c_h(\beta_j)\ c_h(B_k)\ t_h(\beta_k)$, (e)

$1 = -c_h(B_j)\ c_h(\beta_j)\ t_h(B_k)\ c_h(\beta_k)$, (f)

$1 = c_h(B_j)\ c_h(\beta_j)\ t_h(B_k)\ t_h(\beta_k)$ or (g)

$1 = c_h(B_j)\ c_h(\beta_j)\ c_h(B_k)\ c_h(\beta_k)$. (h)

Rewriting equations (A-7) in terms of $\beta_3$ and $B_1$ yields:

$$c_h(\beta_3) = k_c\ c_h(B_1) \text{ or} \tag{A-8}$$

$$c_h(\beta_3) = k_t\ t_h(B_1)\ , \tag{A-9}$$

where the constants $k_c$ and $k_t$ can take on one of the values:

$k_c = t_h(\beta_1)\ c_h(B_3)\ |\ -t_h(\beta_1)\ t_h(B_3)\ |\ -c_h(\beta_1)\ c_h(B_3)\ |\ c_h(\beta_1)\ t_h(B_3)$ or

$k_t = c_h(\beta_1)\ c_h(B_3)\ |\ -c_h(\beta_1)\ t_h(B_3)\ |\ -t_h(\beta_1)\ c_h(B_3)\ |\ t_h(\beta_1)\ t_h(B_3)$.

Converting equations (A-2) and (A-3) with half angle identities and combining with equations (A-8) and (A-9) yields the quadratic equation:

$$A\ c_h^2(\beta_3) + B\ c_h(\beta_3) + C = 0. \tag{A-10}$$

In equation (A-10) the constants A, B and C are defined for four different cases:

$A = (k_c - a)$, $B = -2bk_c$, $C = k_c(ak_c - 1)$ from equations (A-2) and (A-8),

$A = (k_c + a)$, $B = 2bk_c$, $C = -k_c(ak_c + 1)$ from equations (A-3) and (A-8),

$A = (k_t + a)$, $B = -2bk_t$, $C = -k_t(ak_t + 1)$ from equations (A-2) and (A-9), and

$A = (k_t - a)$, $B = 2bk_t$, $C = k_t(ak_t - 1)$ from equations (A-3) and (A-9).

With $\beta_3$ taken from the solution of equation (A-10) the remaining faces angles are known directly from the OAS and OAE assignments and edge lengths from law of sines; the octahedron is completely parameterized.

## Appendix B.

This appendix describes an approach that can be applied to complete the characterization of suspensions of sub-types III-OAE and III-OAS from the parameter set $\{l_1,\alpha_1,\beta_1,\alpha_2,\beta_2,\alpha_3,\alpha_5,\ldots\alpha_{2k+1},\ldots\alpha_{N-3};\ k=1..M-2\}$ when N=2M>4.To accomplish this we generalize relationships that exist between the dihedral angles of vertexes of known flexible suspensions. For Bricard octahedra of the third type all vertexes that are opposite of one another have this relationship, dihedral angles are conjugate or equal in value, and specifically the vertex pair ($\mathbf{v_1},\mathbf{v_3}$) as discussed in Appendix A. In the dodecahedral suspension described in [5] similar relationships are observed between vertex pairs ($\mathbf{v}_1$,$\mathbf{v}_3$), ($\mathbf{v}_2$,$\mathbf{v}_5$) and ($\mathbf{v}_4$,$\mathbf{v}_6$). For N>6 we consider the following vertex pairs to be the generalization of interest:

($\mathbf{v}_1$,$\mathbf{v}_3$),($\mathbf{v}_2$,$\mathbf{v}_5$),($\mathbf{v}_4$,$\mathbf{v}_7$),….($\mathbf{v}_{2k}$,$\mathbf{v}_{2k+3}$),…($\mathbf{v}_{N-4}$,$\mathbf{v}_{N-1}$),($\mathbf{v}_{N-2}$,$\mathbf{v}_N$), k=1..M-2.

For these vertex pairs we develop a recursive approach to the parameterization problem. As an initial step in the recursion the parameter set $\{l_1,\alpha_1,\beta_1,\alpha_2,\beta_2\}$ is used to fully define the four faces $\mathbf{uv}_2\mathbf{v}_1$, $\mathbf{uv}_3\mathbf{v}_2$, $\mathbf{wv}_1\mathbf{v}_2$ and $\mathbf{wv}_2\mathbf{v}_3$ of the suspension as described in Appendix A. The recursion consists of treating $\alpha_{J-2}$ (for J=5,7…N-1) as a parameter and successively computing groups of four adjacent faces by the following computations. The geometry is described in Fig. B-1 where the index substitution J=2k+3 is used to simplify the notation.

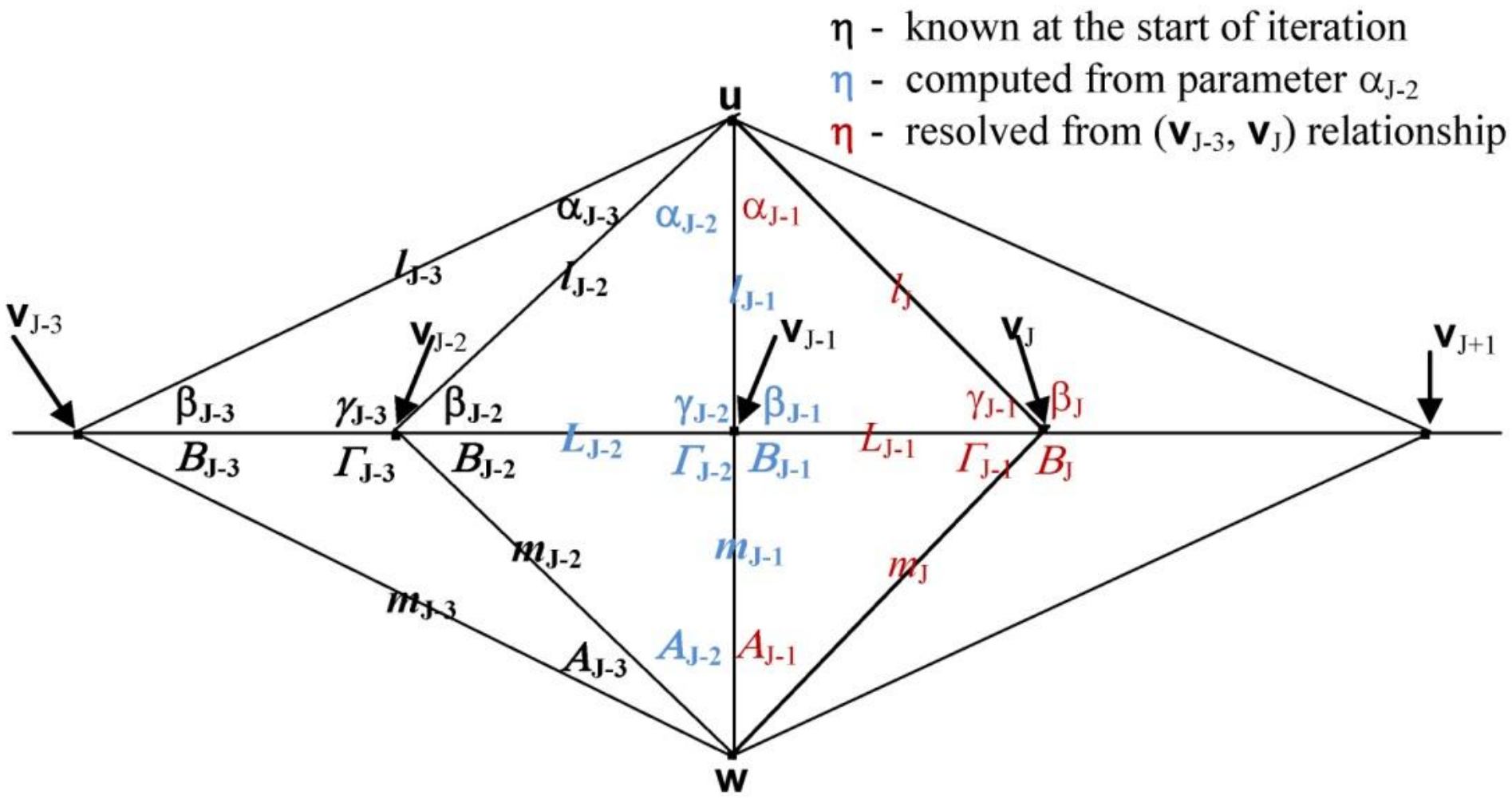

**Figure B-1.** Geometry for Recursive Parameterization

At the start of each step of the recursion, faces $\mathbf{uv}_{J-2}\mathbf{v}_{J-3}$ and $\mathbf{wv}_{J-3}\mathbf{v}_{J-2}$ are known from the previous step. The adjacent faces $\mathbf{uv}_{J-1}\mathbf{v}_{J-2}$ and $\mathbf{wv}_{J-2}\mathbf{v}_{J-1}$ are determined by straightforward geometric computations from the definition of the parameter $\alpha_{J-2}$. The two adjacent faces $\mathbf{uv}_J\mathbf{v}_{J-1}$ and $\mathbf{wv}_{J-1}\mathbf{v}_J$ are defined by using the relationship that exists between the vertex pair ($\mathbf{v}_{J-3}$,$\mathbf{v}_J$) and law of sines equations applied to these faces. Effectively the face angles at vertex $\mathbf{v}_J$ are determined from the face angles at $\mathbf{v}_{J-3}$.

The equivalence relationships between vertex pairs from Appendix A, equations (A-7), are reformulated with j=J and k=J-3 here in terms of $B_J$ and $\beta_J$ as:

$$c_h(B_J) = k_c\, c_h(\beta_J) \text{ or} \tag{B-1}$$

$$c_h(B_J) = k_t\, t_h(\beta_J)\ , \tag{B-2}$$

where the constants $k_c$ and $k_t$ take on one of four values:

$$k_c = c_h(B_{J-3})\, t_h(\beta_{J-3}) \mid t_h(B_{J-3})\, c_h(\beta_{J-3}) \mid -t_h(B_{J-3})\, t_h(\beta_{J-3}) \mid -c_h(B_{J-3})\, c_h(\beta_{J-3}),$$

$$k_t = -t_h(B_{J-3})\, c_h(\beta_{J-3}) \mid -c_h(B_{J-3})\, t_h(\beta_{J-3}) \mid c_h(B_{J-3})\, c_h(\beta_{J-3}) \mid t_h(B_{J-3})\, t_h(\beta_{J-3}).$$

Applying the law of sines to faces **$uv_Jv_{J-1}$** and **$wv_{J-1}v_J$** yields:

$$\frac{l_{J-1}\sin(\beta_{J-1}+\gamma_{J-1})}{\sin\gamma_{J-1}} = \frac{m_{J-1}\sin(B_{J-1}+\Gamma_{J-1})}{\sin\Gamma_{J-1}}.$$

When vertex **$v_J$** is OAE $\gamma_{J-1} = B_J$; and $\Gamma_{J-1} = \beta_J$; this equation reduces to:

$$\cot B_J = a\cot\beta_J + b, \qquad \text{(B-3)}$$

and when vertex **$v_J$** is OAS $\gamma_{J-1} = \pi - B_J$; and $\Gamma_{J-1} = \pi-\beta_J$:

$$\cot B_J = a\cot\beta_J - b, \qquad \text{(B-4)}$$

where

$$a = \frac{m_{J-1}\sin B_{J-1}}{l_{J-1}\sin\beta_{J-1}}$$

$$b = \frac{m_{J-1}\cos B_{J-1} - l_{J-1}\cos\beta_{J-1}}{l_{J-1}\sin\beta_{J-1}}.$$

Combining equations (B-1) through (B-4) yields the quadratic equation:

$$A\,c_h^2(B_J) + B\,c_h(B_J) + C = 0, \qquad \text{(B-5)}$$

where the constants A, B and C are defined for four different cases:

$A = (k_c - a)$, $B = -2bk_c$, $C = k_c(ak_c - 1)$ from equations (B-1) and (B-3),

$A = (k_c - a)$, $B = 2bk_c$, $C = k_c(ak_c - 1)$ from equations (B-1) and (B-4),

$A = (k_t + a)$, $B = -2bk_t$, $C = -k_t(ak_t + 1)$ from equations (B-2) and (B-3), and

$A = (k_t + a)$, $B = 2bk_t$, $C = -k_t(ak_t + 1)$ from equations (B-2) and (B-4).

Faces **$uv_Jv_{J-1}$** and **$wv_{J-1}v_J$** are fully defined using $B_J$ from the solution of equation (B-5) and $\beta_J$ from equations (B-3) or (B-4).

The above recursion defines four faces from a single parameter and existing lower indexed faces. From a practical standpoint the multiplicity of choices of possible solutions of Eq. (B-5) at each stage can present a difficult computational problem. To constrain the number of choices at each step of the recursion we add the notion of a partial construction at each stage in which the derived face angles and edge lengths satisfy not only the described recursion but also satisfy conditions that are imposed by Eq. (A-4) by the construction at a specific and convenient value of the variable of flexion, such as $\varepsilon_1 = \pi/2$. By using this approach, derived parameters that do not lead to realizable flexible suspensions are eliminated.

## Appendix C.

This appendix describes useful relationships that exist between dihedral angles and face angles of a cap of index N. These relationships can be developed from a consideration of such a cap that has unit length edges at the apical vertexes and is oriented as shown in Fig. C-1with the apical vertex located at the origin of coordinates. From inspection of this figure it is apparent that the position of the k-th vertex can be achieved by either a series of rotations starting with the first vertex and progressing through the second, third etc. vertexes or alternately by starting with the N-th and progressing through the N-th, (N-1)-th, (N-2)-th etc. vertexes.

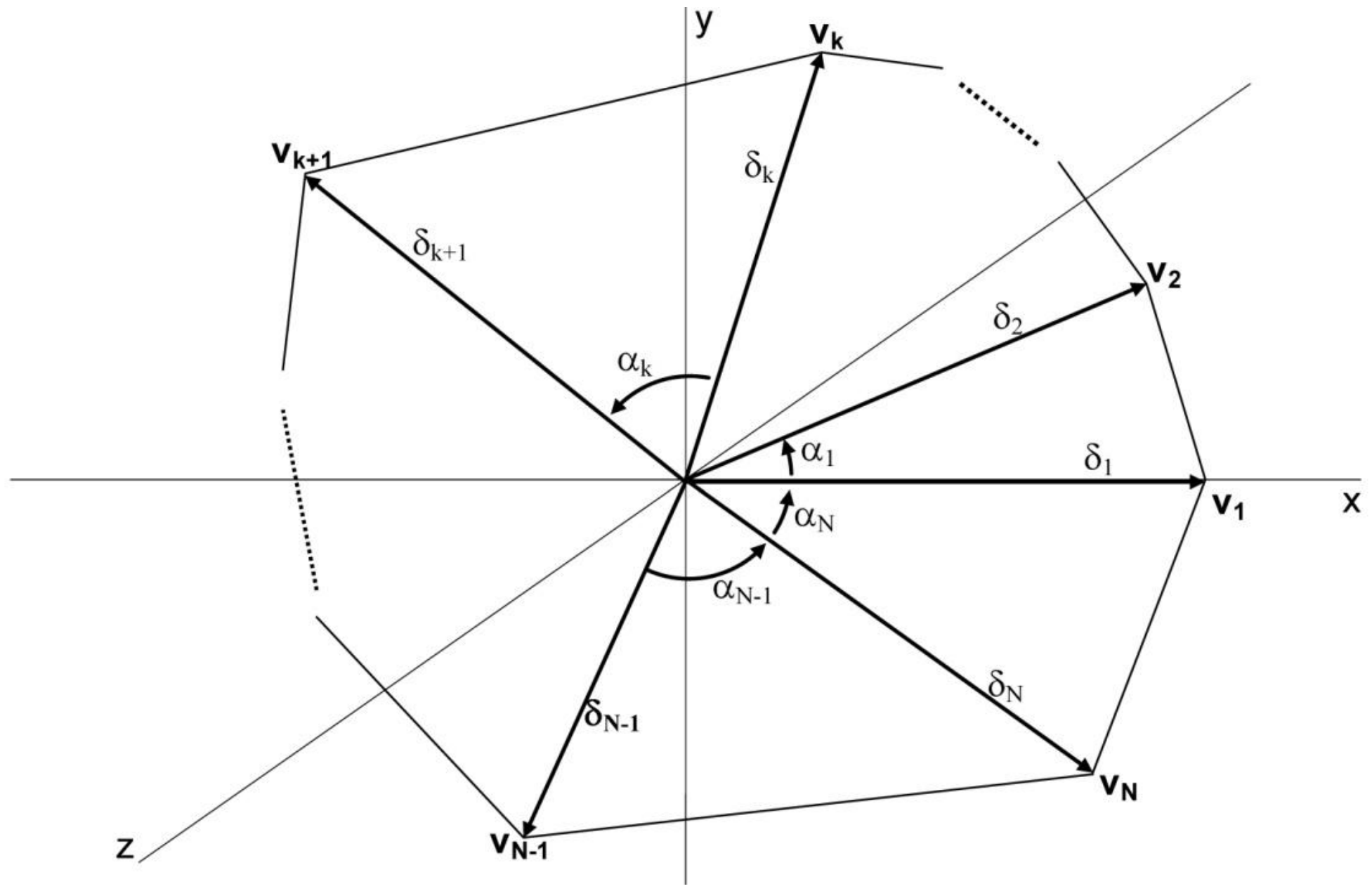

**Figure C-1.** Dihedral/Face Angle Relationships

The first rotation sequence is given by:

$$\mathbf{v_1} = \begin{pmatrix} 1 \\ 0 \\ 0 \end{pmatrix},$$

$$\mathbf{v_2} = A(\alpha_1)\,\mathbf{v_1},$$
$$\mathbf{v_3} = A(\alpha_1)\,\Phi(\pi-\delta_2)\,A(\alpha_2)\,\mathbf{v_1},$$
$$\mathbf{v_4} = A(\alpha_1)\,\Phi(\pi-\delta_2)\,A(\alpha_2)\,\Phi(\pi-\delta_3)\,A(\alpha_3)\,\mathbf{v_1},$$
$$\cdot$$
$$\cdot$$
$$\mathbf{v_k} = A(\alpha_1)\prod_{i=2..k-1}[\Phi\,(\pi-\delta_i)\,A(\alpha_i)]\,\mathbf{v_1} \text{ for } k \leq N, \qquad \text{(C-1)}$$

while the second is given by:

$$\mathbf{v_N} = \Phi(\delta_1-\pi)\, A(-\alpha_N)\, \mathbf{v_1}.$$
$$\mathbf{v_{N-1}} = \Phi(\delta_1-\pi)\, A(-\alpha_N)\, \Phi(\delta_N-\pi)\, A(-\alpha_{N-1})\, \mathbf{v_1}.$$
$$\cdot$$
$$\cdot$$
$$\mathbf{v_k} = \Phi(\delta_1-\pi)\, A(-\alpha_N) \prod_{i=N..k+1} [\Phi(\delta_i-\pi)\, A(-\alpha_{i-1})]\, \mathbf{v_1} \text{ for } k \geq 1. \quad \text{(C-2)}$$

Here $A(\alpha)$ and $\Phi(\delta)$ are rotation matrixes defined by

$$A(\alpha) = \begin{pmatrix} \cos\alpha & -\sin\alpha & 0 \\ \sin\alpha & \cos\alpha & 0 \\ 0 & 0 & 1 \end{pmatrix},$$

and

$$\Phi(\delta) = \begin{pmatrix} 1 & 0 & 0 \\ 0 & \cos\delta & -\sin\delta \\ 0 & \sin\delta & \cos\delta \end{pmatrix}.$$

Since there are N-3 independent dihedral variables associated with this cap it is of interest to be able to determine the three dependent dihedrals from these equations under the assumption that the face angles$\{\alpha_i,\ i=1..N\}$are completely defined and that of the dihedrals $\{\delta_i,\ i=1..N\}$, the three $\delta_{N-2}$, $\delta_{N-1}$ and $\delta_N$ are unknown.

Equating the two forms of the equations (C-1) and (C-2) for the vertex $\mathbf{v_{N-2}}$ yields:

$$\Phi(\delta_1-\pi)A(-\alpha_N) \prod_{i=N..N-1} [\Phi(\delta_i-\pi)A(-\alpha_{i-1})]\mathbf{v_1} = A(\alpha_1) \prod_{i=2..N-3} [\Phi(\pi-\delta_i)\, A(\alpha_i)]\mathbf{v_1}.$$

After multiplication on the left by the matrix product $A(\alpha_{N-1})\Phi(\pi-\delta_N)\, A(\alpha_N)\, \Phi(\pi-\delta_1)$ this reduces to:

$$\Phi(\delta_{N-1}-\pi)A(-\alpha_{N-2})\,\mathbf{v_1} = A(\alpha_{N-1})\Phi(\pi-\delta_N)\,\mathbf{M}\,\mathbf{v_1} \quad \text{(C-3)}$$

where the matrix $\mathbf{M}$ is defined as follows:

$$\mathbf{M} = A(\alpha_N) \prod_{i=1..N-3} [\Phi(\pi-\delta_i)\, A(\alpha_i)].$$

The three components of equation (C-3) evaluate to the equations:

$$\begin{aligned} \cos\alpha_{N-2} &= M_{11} \cos\alpha_{N-1} + \sin\alpha_{N-1}\, (M_{12} \cos\delta_N + M_{13} \sin\delta_N), \\ \sin\alpha_{N-2} \cos\delta_{N-1} &= M_{11} \sin\alpha_{N-1} - \cos\alpha_{N-1}\, (M_{12} \cos\delta_N + M_{13} \sin\delta_N) \text{ and} \\ \sin\alpha_{N-2} \sin\delta_{N-1} &= M_{12} \sin\delta_N - M_{13} \cos\delta_N. \end{aligned} \quad \text{(C-4)}$$

These equations can be solved for the dihedral angles $\delta_{N-1}$ and $\delta_N$. The remaining dihedral angle $\delta_{N-2}$ can be found by similar manipulations of the two forms of the vertex $\mathbf{v_{N-1}}$ or of the the vertex $\mathbf{v_N}$. Using the former leads to:

$$A(\alpha_1) \prod_{i=2..N-2} [\Phi\, (\pi-\delta_i)\, A(\alpha_i)]\, \mathbf{v_1} = \Phi(\delta_1-\pi)\, A(-\alpha_N)\, \Phi(\delta_N-\pi)\, A(-\alpha_{N-1})\, \mathbf{v_1}.$$

Multiplication on the left by $A(\alpha_N)\, \Phi(\delta_1-\pi)$ leads to:

$$\mathbf{M}\,\Phi(\pi-\delta_{N-2})\, A(\alpha_{N-2})\mathbf{v_1} = \Phi(\delta_N-\pi)\, A(-\alpha_{N-1})\, \mathbf{v_1}$$

which can be readily solved for $\delta_{N-2}$ from the y and z components ($\mathbf{P} = \mathbf{M}^{-1}$):

$$\sin\alpha_{N-2} \cos\delta_{N-2} = -P_{21} \cos\alpha_{N-1} - \sin\alpha_{N-1} (P_{22} \cos\delta_N + P_{23} \sin\delta_N) \text{ and}$$
$$\sin\alpha_{N-2} \sin\delta_{N-2} = P_{31} \cos\alpha_{N-1} + \sin\alpha_{N-1} (P_{32} \cos\delta_N + P_{33} \sin\delta_N).$$

For the case N=4, the x-component of the two forms of $\mathbf{v_4}$, when equated and after some re-arrangement of terms, yields the equation:

$$\begin{aligned}&\sin\alpha_1 \cos\alpha_2 \sin\alpha_3 \cos\delta_2 \cos\delta_3 - \sin\alpha_1 \sin\alpha_3 \sin\delta_2 \sin\delta_3 - \\ &\quad \sin\alpha_1 \sin\alpha_2 \cos\alpha_3 \cos\delta_2 - \cos\alpha_1 \sin\alpha_2 \sin\alpha_3 \cos\delta_3 + \\ &\qquad \cos\alpha_4 - \cos\alpha_1 \cos\alpha_2 \cos\alpha_3 = 0.\end{aligned} \quad \text{(C-5)}$$

By using the following association of notation: $\alpha_1 \Leftrightarrow \delta$, $\alpha_2 \Leftrightarrow \alpha$, $\alpha_3 \Leftrightarrow \beta$, $\alpha_4 \Leftrightarrow \gamma$, $\delta_2 \Leftrightarrow \psi$ and $\delta_3 \Leftrightarrow \phi$; the above equation is identical to the equation derived in [1] to describe the deformation of a cap with four faces. There it is transformed by the transformations $t = \tan(\phi/2)$ and $u = \tan(\psi/2)$ into the following "quadratic" tetrahedral equation [1, Eq. 1] that forms the basis for the analysis contained therein:

$$At^2u^2 + Bt^2 + 2Ctu + Du^2 + E = 0.$$

Here the coefficients are defined by

$$\begin{aligned}A &= \cos(\gamma) - \cos(\alpha+\beta+\delta),\\ B &= \cos(\gamma) - \cos(\alpha+\beta-\delta),\\ C &= -2\sin(\beta)\sin(\delta),\\ D &= \cos(\gamma) - \cos(\alpha-\beta+\delta) \text{ and}\\ E &= \cos(\gamma) - \cos(\alpha-\beta-\delta).\end{aligned}$$

Similarly, the x-component of the two forms of $\mathbf{v_3}$, when equated yields the equation

$$\cos\alpha_1 \cos\alpha_2 + \sin\alpha_1 \sin\alpha_2 \cos\delta_2 = \cos\alpha_3 \cos\alpha_4 + \sin\alpha_3 \sin\alpha_4 \cos\delta_4. \quad \text{(C-6)}$$

This is the "linear" equation that relates opposite dihedral angles at a vertex of index 4 and is also used in [1, Sec. 4] but is not explicitly identified there.

## Appendix D.

This appendix contains details of example flexible suspensions that illustrate the generalizations contained in Secs. 4 through 9 and that are constructed using the approach described in Sec. 10. Also included are images of selected examples along with a description of the primary computational algorithm used to compute dihedral angles and identify specific constructions. Additionally, several ancillary computations that are used to validate the flexibility of the example suspensions are identified. The existance of rigid isomers is mentioned briefly.

**Parameter Sets:** Parameter sets for values of N=4,6,8,10,12,14,16 are provided for all five sub-types I-OEE, II-AEE, II-OEE, III-OAE and III-OAS. All values in these parameter sets are integer values to simplify the data tables and this fact is of no particular significance as the constructions do not require this attribute. All values can be varied over some non-integer interval individually, or collectively, to sucessfully construct similar example suspensions. Also the occurrence of parameters having the same value is generally not significant. In cases where it is relevant, say to create symmetry, a different notation is used in the tables; eg. $L_1$=$L_3$=5 rather than $L_1$=5 and $L_3$=5.

For most cases a single parameter set represents more than one flexible suspension. For example, a parameter set for sub-type I-OEE with N=10 represents 16 ($=2^4=2^{M-1}$) flexible suspensions, each of which is constructed with a different folding. Similarly for sub-types II-AEE or II-OEE. On the other hand, for sub-types III-OAE and III-OAS each parameter set may represent several flexible suspensions that have different computed face angles and edge lengths but each of which has only one folding that is flexible. Consequently in the subsequent discussion of specific examples the suspension of interest will be identified by the associated parameter set and a Dihedral Identifier (DI) that represents the specific folding that is used to create the suspension. The significance of this identifier is described in the following paragraphs.

**Dihedral Identifier:** As discussed in Sec. 10 the primary technique that is used to construct suspensions is the computation of the dihedral angles at vertex $\mathbf{v_k}$ based upon dihedral angle $\varepsilon_{k-1}$ at vertex $\mathbf{v_{k-1}}$ given that all face angles at vertex $\mathbf{v_k}$ are known. This requires the solution of equation (C-5) which is quadratic and thereby introduces the choice of sign in the solution of the quadratic equation as a construction option. The DI is the catenation of one's (plus solutions) and zero's (negative solutions) into a composite integer value in which the value for the k-th solution is placed in the k-th bit of the DI. The quadratic equation that a particular bit applies to is created from equation (C-5) as applied to vertex $\mathbf{v_k}$ by using the following association of face angle and dihedral angle notation:

$$\alpha_1 \Leftrightarrow \beta_k,\ \alpha_2 \Leftrightarrow \gamma_{k-1},\ \alpha_3 \Leftrightarrow \mathit{\Gamma}_{k-1},\ \alpha_4 \Leftrightarrow \mathit{B}_k,\ \delta_2 \Leftrightarrow \delta_k \text{ and } \delta_3 \Leftrightarrow \varepsilon_{k-1}.$$

With this association equation (C-5) may be written as

$$A \cos\delta_k \cos\varepsilon_{k-1} + B \sin\delta_k \sin\varepsilon_{k-1} + C \cos\delta_k + D \cos\varepsilon_{k-1} + E = 0,$$

where the coefficients are defined by

$$A = \sin\beta_k \cos\gamma_{k-1} \sin\mathit{\Gamma}_{k-1},$$
$$B = -\sin\beta_k \sin\mathit{\Gamma}_{k-1},$$
$$C = -\sin\beta_k \sin\gamma_{k-1} \cos\mathit{\Gamma}_{k-1},$$
$$D = -\cos\beta_k \sin\gamma_{k-1} \sin\mathit{\Gamma}_{k-1} \text{ and}$$
$$E = \cos\mathit{B}_k - \cos\beta_k \cos\gamma_{k-1} \cos\mathit{\Gamma}_{k-1}.$$

By treating $\varepsilon_{k-1}$ as known and introducing half angle equations for $\delta_k$ this reduces to:

$$a\, t_h^2(\delta_k) + b\, t_h(\delta_k) + c = 0, \qquad \text{(D-1)}$$

where the function $t_h(\delta)$ is defined in equations (A-5) and the coefficients a, b and c are defined by:

$$a = (E-C) + (D-A) \cos\varepsilon_{k-1},$$
$$b = 2B \sin\varepsilon_{k-1} \text{ and}$$
$$c = (E+C) + (D+A) \cos\varepsilon_{k-1}.$$

The sign that is used to solve the quadratic equation (D-1) is retained in the appropriate bit of the DI. The dihedral angle $\varepsilon_k$ at vertex $\mathbf{v_k}$ is determined by a straight forward application of equation (C-6) using the dihedral $\varepsilon_{k-1}$ that was used in the solution of equation (D-1). The DI's presented in the subsequent discussion represent the solutions for all N non-apical vertexes.

**Dihedral Angles:** Dihedral angles computed from equation (D-1) are adequate for use in positioning of vertexes but require some adjustments when considered over the full range of flexion and when used in the computation of some related quantities of interest. Without adjustments these dihedral angles will not be continuous over the full range of flexion. Generally, dihedral angles of suspensions are considered to be multivalued functions [5, pg. 20] of the form $\delta=\delta_0+2\pi k$ for some integer value k and $\delta_0$ is the "normal" value in the range $[0,2\pi]$. We deal with this problem as follows.

The variable of flexion for all suspension sub-types is the dihedral angle $\varepsilon_1$. For suspensions of sub-type I-OEE, II-AEE and II-OEE the range (of flexibility) of $\varepsilon_1$ is found to be in one of three forms: two intervals, $[\varepsilon_m, \varepsilon_x]$ and $[\varepsilon_M, \varepsilon_X]$ where $2\pi>\varepsilon_X>\varepsilon_M>\pi>\varepsilon_x>\varepsilon_m>0$ or two intervals in which $\varepsilon_m=0$ and $\varepsilon_X=2\pi$. In the latter case $[\varepsilon_M, \varepsilon_x+2\pi]$ is taken as the interval of interest. The third form has a single interval $[\varepsilon_m, \varepsilon_x]$ with $2\pi\geq\varepsilon_x>\pi>\varepsilon_m\geq 0$. Since the range of flexion for suspensions of sub-type III-OAE or III-OAS is $[0, 2\pi]$ no adjustments to $\varepsilon_1$ are required for these sub-types. All other dihedrals are treated as continuous non-negative functions of $\varepsilon_1$ over the range, or ranges, of flexion. This necessitates adjustments by $\pm 2\pi$ whenever the value of a specific dihedral passes from a value near 0 to a value near $2\pi$ and vice versa.

The adjusted dihedral angle values are used in a variety of computations that validate the construction of the suspension and its flexibility. These computations consist of 1) the derivative term, equation (10), is evaluated for the edge length $L_N=|\mathbf{v_N v_1}|$, 2) the oriented volume of the suspension and 3) the total mean curvature [9]; the latter of which is defined [10] by terms: (edge length)*( π– dihedral) summed over all edges of the suspension **S**. Using the notation of this paper:

$$2C(\mathbf{S}) = \sum_{k=1..N} [l_k(\pi - \delta_k) + m_k(\pi - \Delta_k) + L_k(\pi - \varepsilon_k)].$$

Additionally the relationships between the dihedral angles of the suspensions given by equations (2), (4), (6) and (8) effectively define computational constants that are convenient to evaluate. As a consequence of these relationships the sum of all dihedral angles in a flexible suspension is constant under flexion.

**Solid Angles:** There exist a number of relationships between the solid angles (computed with Girard's Theorem) of vertex pairs that are constant during flexion. Typically the solid angles of each of the vertex pairs add together to $=4\pi$. The vertex pairs of interest are sub-type dependent as follows: sub-types I-OEE and I-OEE have vertex pairs $(\mathbf{v_k}, \mathbf{v_{k+M}})$ for k=1..M, sub-type II-OAE has $(\mathbf{v_1}, \mathbf{v_{M+1}})$ and $(\mathbf{v_k}, \mathbf{v_{N-k+1}})$ for k=2..M, and III-OAE and III-OAS have the pairs $\{(\mathbf{v_1},\mathbf{v_3}), (\mathbf{v_{N-2}},\mathbf{v_N}), (\mathbf{v_{2k}},\mathbf{v_{2k+3}})$ , k=1..M-2}. In all suspensions the apical vertexes **u** and **w** are paired together. As a consequence the sum of all solid angles is constant under flexion.

Although this relationship between the solid angles of the relevant vertex pairs is observed in all examples of all sub-types it is not known if this relationship constitutes a necessary and sufficient condition for flexibility of a suspension.

Further in all but the I-OEE sub-type there are vertex pairs in which the individual vertexes have a constant solid angle, typically $=2\pi$. In the sub-type II-AEE the pair $(\mathbf{v_1}, \mathbf{v_{M+1}})$ has this property while in sub-type II-OEE the pair of interest is (**u**, **w**). The property is also observed in sub-types III-OAE and III-OAS but does not seem to have such an easy quantification.

**Example I-OEE, II-AEE, II-OEE Suspensions:** Tables D-I, D-II and D-III contain parameter sets for the sub-types I-OEE, II-AEE, II-OEE respectively.

Parameter sets are readily created for these three sub-types. Any parameters that satisfy the constraints found in equations (1), (3) and (5) and for which the triangle inequalities are satisfied can be used to create flexible suspensions. The primary computational challenge associated with these types is the determination of the range of values of $\varepsilon_1$ for which the suspension is flexible. For all examples this is accomplished by repeated attempts at construction with $\varepsilon_1$ taken over the interval $[0, 2\pi]$. For all examples the flexibility is evaluated at $1^o$ increments and minimum and maximum values resolved iteratively to $\sim(10^{-4})^o$.

Parameters that define dipyramids satisfy the constraints of sub-types I-OEE, II-AEE, II-OEE and are observed to have flexible foldings that have the interesting feature that all non-apical vertexes are co-planar under flexion. However all that we have observed are also degenerate in the sense that some non-apical vertexes are always coincident with one another.

| # | N | Parameters |
|---|---|---|
| **1** | **4** | $l_1=10,\ l_2=11,\ l_3=12,\ l_4=13;\ L_1=8,\ L_2=9.$ |
| **2** | **6** | $l_1=10,\ l_2=12,\ l_3=11,\ l_4=12.\ l_5=11,\ l_6=13;\ L_1=3,\ L_2=4,\ L_3=5.$ |
| **3** | **8** | $l_1=10,\ l_2=12,\ l_3=11,\ l_4=13,\ l_5=14,\ l_6=16,\ l_7=15,\ l_8=14;\ L_1=4,\ L_2=3,\ L_3=5,\ L_4=6.$ |
| **4** | **10** | $l_1=10,\ l_2=13,\ l_3=15,\ l_4=13,\ l_5=16,\ l_6=12,\ l_7=15,\ l_8=14,\ l_9=12,\ l_{10}=13;\ L_1=10,\ L_2=11,\ L_3=12,\ L_4=14,\ L_5=13.$ |
| **5** | **12** | $l_1=10,\ l_2=13,\ l_3=15,\ l_4=13,\ l_5=16,\ l_6=12,\ l_7=15,\ l_8=14,\ l_9=12,\ l_{10}=13,\ l_{11}=14,\ l_{12}=11;$<br>$L_1=10,\ L_2=11,\ L_3=12,\ L_4=14,\ L_5=13,\ L_6=12.$ |
| **6** | **14** | $l_1=10,\ l_2=13,\ l_3=15,\ l_4=13,\ l_5=16,\ l_6=12,\ l_7=15,\ l_8=14,\ l_9=12,\ l_{10}=13,\ l_{11}=14,\ l_{12}=11,\ l_{13}=13,\ l_{14}=12;$<br>$L_1=10,\ L_2=15,\ L_3=12,\ L_4=14,\ L_5=13,\ L_6=12,\ L_7=11.$ |
| **7** | **16** | $l_1=10,\ l_2=11,\ l_3=12,\ l_4=10,\ l_5=10,\ l_6=12,\ l_7=11,\ l_8=12,\ l_9=10,\ l_{10}=12,\ l_{11}=11,\ l_{12}=12,\ l_{13}=10,\ l_{14}=12,\ l_{15}=11,\ l_{16}=12;$<br>$L_1=10,\ L_2=11,\ L_3=12,\ L_4=10,\ L_5=12,\ L_6=10,\ L_7=12,\ L_8=11.$ |

**Table D-I.** I-OEE Suspension Example Parameter Sets

| # | N | Parameters |
|---|---|---|
| **1** | **4** | $l_1=10,\ l_2=11,\ l_3=12,\ l_4=13;\ L_1=5,\ L_2=4.$ |
| **2** | **6** | $l_1=10,\ l_2=14,\ l_3=15,\ l_4=12.\ l_5=13,\ l_6=11;\ L_1=5,\ L_2=4,\ L_3=6.$ |
| **3** | **8** | $l_1=10,\ l_2=11,\ l_3=12,\ l_4=13,\ l_5=12,\ l_6=11,\ l_7=14,\ l_8=13;\ L_1=5,\ L_2=6,\ L_3=4,\ L_4=3.$ |
| **4** | **10** | $l_1=10,\ l_2=13,\ l_3=14,\ l_4=14,\ l_5=12,\ l_6=11,\ l_7=13,\ l_8=13,\ l_9=17,\ l_{10}=16;\ L_1=10,\ L_2=9,\ L_3=8,\ L_4=7,\ L_5=7.$ |
| **5** | **12** | $l_1=10,\ l_2=12,\ l_3=11,\ l_4=14,\ l_5=13,\ l_6=16,\ l_7=15,\ l_8=18,\ l_9=17,\ l_{10}=20,\ l_{11}=19,\ l_{12}=22;$<br>$L_1=15,\ L_2=16,\ L_3=17,\ L_4=18,\ L_5=19,\ L_6=20.$ |
| **6** | **14** | $l_1=10,\ l_2=12,\ l_3=11,\ l_4=14,\ l_5=13,\ l_6=16,\ l_7=15,\ l_8=18,\ l_9=17,\ l_{10}=20,\ l_{11}=19,\ l_{12}=22,\ l_{13}=21,\ l_{14}=24;$<br>$L_1=16,\ L_2=17,\ L_3=18,\ L_4=19,\ L_5=20\ L_6=21,\ L_7=22.$ |
| **7** | **16** | $l_1=10,\ l_2=15,\ l_3=20,\ l_4=11,\ l_5=16,\ l_6=21,\ l_7=12,\ l_8=17,\ l_9=22,\ l_{10}=13,\ l_{11}=18,\ l_{12}=23,\ l_{13}=14,\ l_{14}=19,\ l_{15}=24,\ l_{16}=17;$<br>$L_1=15,\ L_2=16,\ L_3=17,\ L_4=18,\ L_5=18,\ L_6=17,\ L_7=16,\ L_8=15.$ |

**Table D-II.** II-AEE Suspension Example Parameter Sets

| # | N | Parameters |
|---|---|---|
| **1** | **4** | $l_1=10,\ l_2=13;\ m_1=16,\ m_2=12;\ L_1=8,\ L_2=7.$ |
| **2** | 6 | $l_1=10,\ l_2=11,\ l_3=12;\ m_1=12.\ m_2=13,\ m_3=15;\ L_1=5,\ L_2=3,\ L_3=4.$ |
| **3** | **8** | $l_1=10,\ l_2=12,\ l_3=11,\ l_4=13;\ m_1=14,\ m_2=17,\ m_3=15,\ m_4=13;\ L_1=4,\ L_2=3,\ L_3=5,\ L_4=6.$ |
| **4** | **10** | $l_1=10,\ l_2=13,\ l_3=14,\ l_4=15,\ l_5=12;\ m_1=18,\ m_2=21,\ m_3=19,\ m_4=16,\ m_5=15;\ L_1=9,\ L_2=10,\ L_3=8,\ L_4=6,\ L_5=7.$ |
| **5** | **12** | $l_1=10,\ l_2=13,\ l_3=14,\ l_4=15,\ l_5=12,\ l_6=11;\ m_1=11,\ m_2=14,\ m_3=16,\ m_4=13,\ m_5=13,\ m_6=12;$<br>$L_1=5,\ L_2=4,\ L_3=7,\ L_4=8,\ L_5=3,\ L_6=9.$ |
| **6** | **14** | $l_1=10,\ l_2=13,\ l_3=14,\ l_4=15,\ l_5=12,\ l_6=11,\ l_7=13;\ m_1=11,\ m_2=14,\ m_3=15,\ m_4=16,\ m_5=13,\ m_6=12,\ m_7=14;$<br>$L_1=10,\ L_2=10,\ L_3=8,\ L_4=7,\ L_5=6,\ L_6=9,\ L_7=12.$ |
| **7** | **16** | $l_1=l_2=l_3=l_4=l_5=l_6=l_7=l_8=10;\ m_1=m_2=m_3=m_4=m_5=m_6=m_7=m_8=13;\ L_1=L_3=L_5=L_7=8,\ L_2=L_4=L_6=L_8=10.$ |

**Table D-III.** II-OEE Suspension Example Parameter Sets

In the following pages we present images of several example flexible suspensions. All images are rendered with ray trace software [12] and all suspensions are presented in scenes in front of a magenta tinted checkered wall and above a grey tinted checkered floor. This has been done to provide some perspective when various aspects of the same suspension are viewed from different eye points; eg. a front view versus top view. The dimension of the checkers is 10x10 (unit less) which is the same dimension that is used in all constructions for the length of the edge $l_1$. In all images the outsides of faces are shown in blue; insides in yellow. Self intersections of faces are shown where there is a change of color or shading; edges are delineated by a black line. Motion is depicted by the display of red spheres over the trajectory of a vertex and in some cases by meshes of red lines representing edges at various positions during flexion. Axes and planes of interest are shown in aqua.

An **I-OEE Flexible Suspension** with 28 faces (N=14) that is constructed with $\varepsilon_1=100^o$ from parameter set #6 found in Table D-I is illustrated in Fig. D-1. This suspension has a DI=16129 =$3F01_{16}$. This DI illustrates a general characteristic of the sub-type since the bits satisfy the equation $b_k=1-b_{k+M}$ for k=1..M. This indicates that both solutions of the equation (B-5) are used in the construction. Fig. D-1A shows a "front" view of the suspension with the axis of symmetry oriented vertically in the middle of the image. The linear motion of the apical vertexes is shown along with the more general arc like motion of the symmetric vertexes $\mathbf{v_1}$ and $\mathbf{v_8}$. The insets depict end views of just the respective caps at the apical vertexes and illustrate how the two caps are inverted from one another. They do in fact "fit" into one another as their solid angles are space filling, ie. $\sigma_u+\sigma_w=4\pi$.

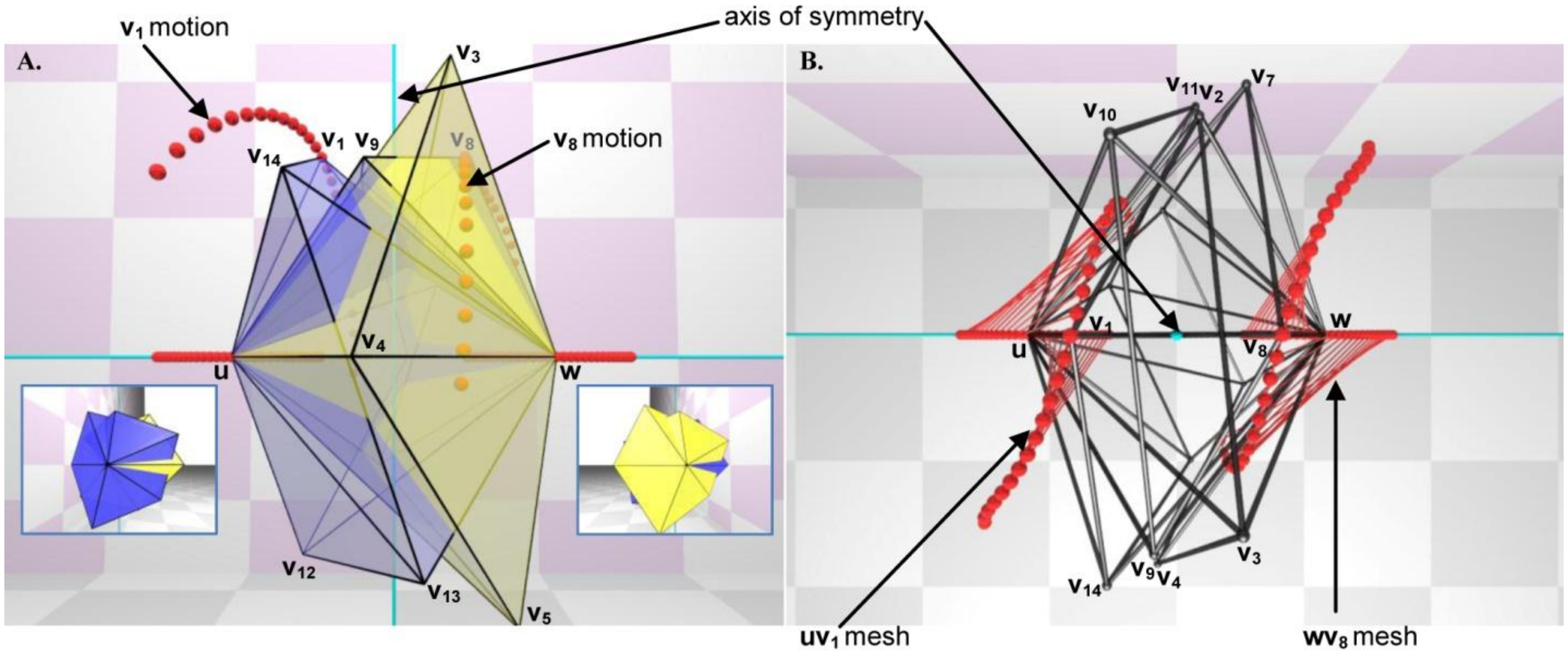

**Figure D-1.** I-OEE Flexible Suspension. **A.** Front View. **B.** Top View.

Fig. D-1B is a "top" view of the I-OEE suspension looking down on the axis of symmetry. Face images have been removed and edges and vertexes shown more distinctly as a framework. The framework of edges and labeled vertexes illustrate the axial symmetry of a specific construction (with $\varepsilon_1=100^o$) while under flexion the axial symmetry is illustrated by the red meshs swept out by the edges, $\mathbf{uv_1}$ and $\mathbf{wv_8}$, that are formed from the symmetrically located vertex pair ($\mathbf{v_1}$,$\mathbf{v_8}$).

An **II-AEE Flexible Suspension** that has 20 faces (N=10) and is constructed with $\varepsilon_1=145^o$ from parameter set #4 found in Table D-II is illustrated in Fig. D-2. This suspension has a DI=301=$12D_{16}$ which illustrates a general characteristic of the sub-type since the bits satisfy the equation $b_k=1-b_{N-k+1}$ for k=1..M. Fig. D-2A shows a "front" view of the II-AEE suspension where the plane of symmetry of this II-AEE suspension appears in the middle of the image and is defined by the motion traces of vertexes $\mathbf{v_1}$ and $\mathbf{v_6}$. The linear motion traces of vertexes **u** and **w** are perpendicular to this plane. Fig. D-1B is a “top” view of the same suspension looking down on the edge of the plane of symmetry. The framework of edges and labeled vertexes illustrate the planar symmetry of a specific construction while under flexion this is illustrated by the red meshs swept out by the edges $\mathbf{uv_2}$ and $\mathbf{wv_{10}}$. The vertexes $\mathbf{v_2}$ and $\mathbf{v_{10}}$ are symmetrically located.

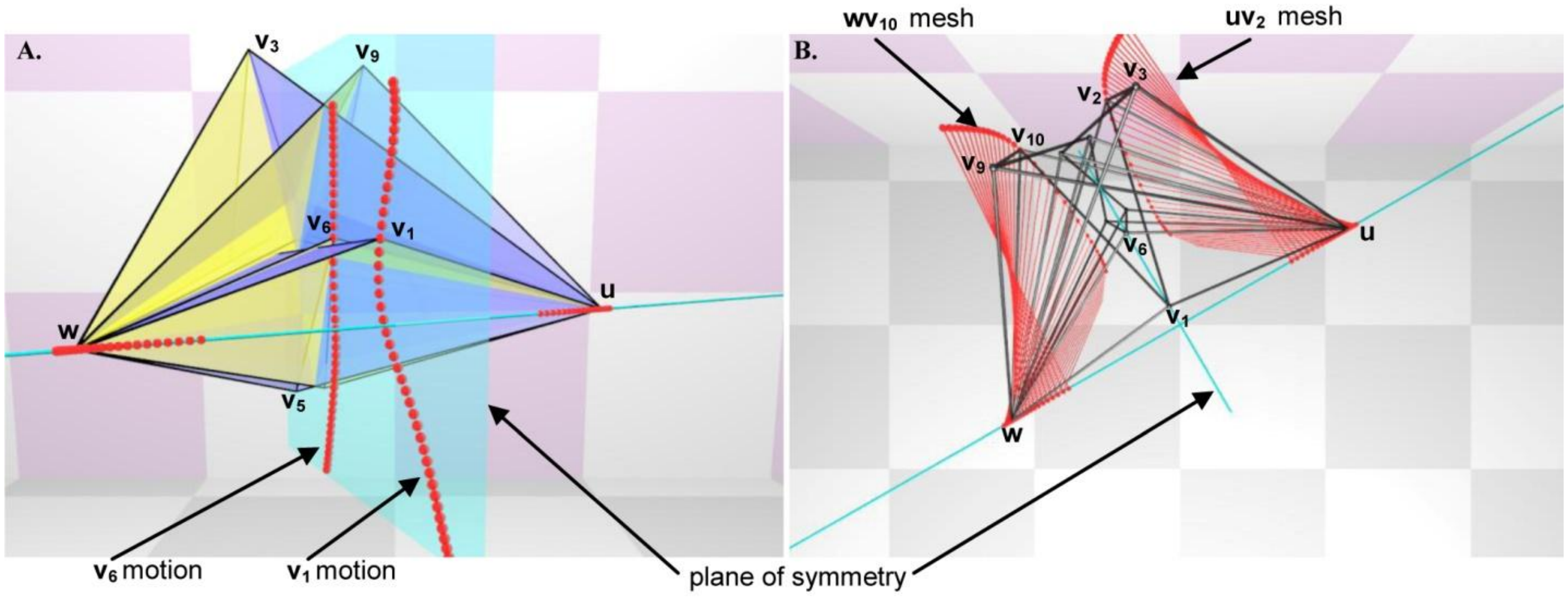

**Figure D-2.** II-AEE Flexible Suspension. **A.** Front View. **B.** Top View.

Fig. D-3 illustrates the construction of an **II-OEE Flexible Suspension** that has 32 faces (N=16) and is based upon parameter set #7 in Table D-III. This suspension is constructed with the variable of flexion $\varepsilon_1$=75$^o$ and DI=41055=A05F$_{16}$. This DI illustrates a general characteristic of the II-OEE sub-type since the bits satisfy the equation $b_k$=1-$b_{k+M}$ for k=1..M.

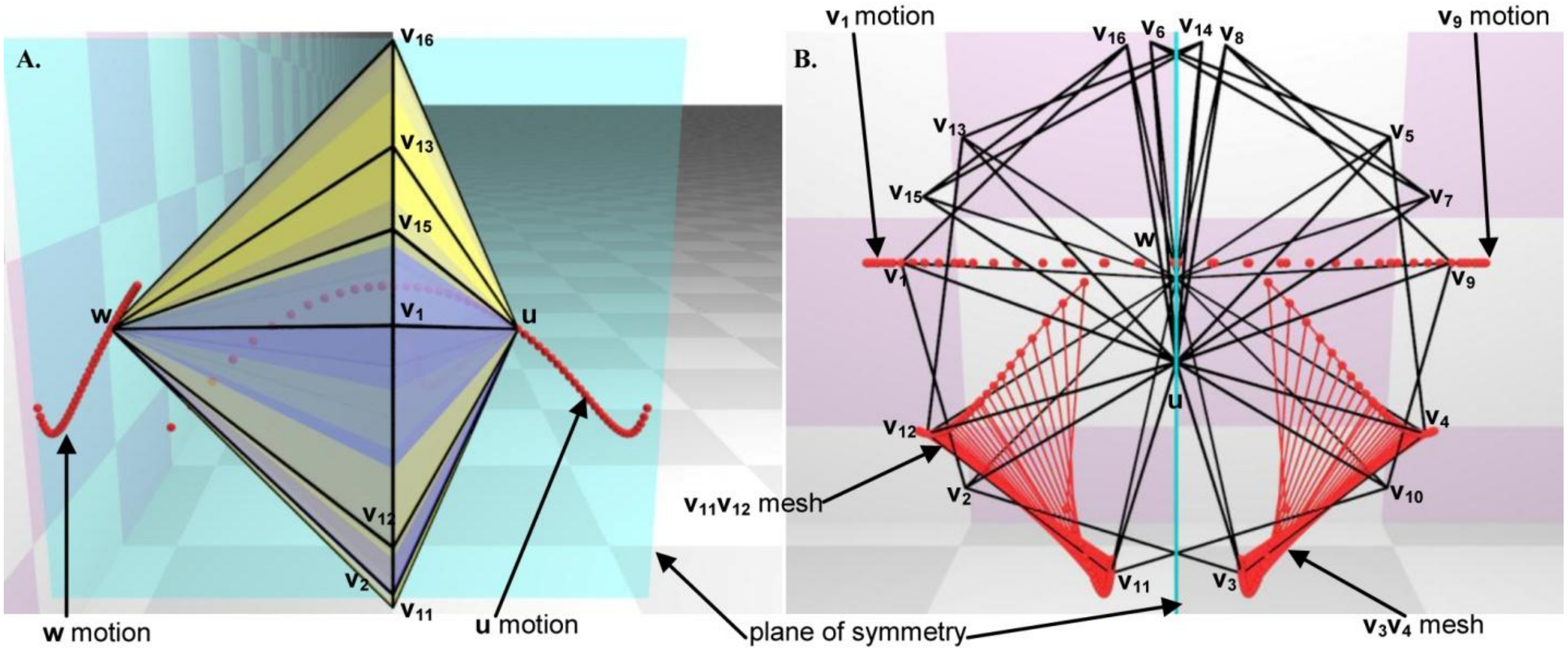

**Figure D-3.** II-OEE Flexible Suspension. **A.** Side View. **B.** Front View.

Fig. D-3A is a "side view" that illustrates the symmetry of the construction and the interesting feature that all non-apical vertexes are co-planar under flexion. The plane of symmetry for the suspension is parallel to the plane of the page and is defined by the motion traces of the apical vertexes **u** and **w**. Fig. D-3B is a "front" view of the II-OEE suspension looking at the edge of the plane of symmetry. Edges and labeled vertexes illustrate the planar symmetry of the specific construction with $\varepsilon_1$=75$^o$. Under flexion this symmetry is illustrated by the red meshs swept out by the edges, $\mathbf{v_3v_4}$ and $\mathbf{v_{11}v_{12}}$, that are formed from symmetrically located vertex pairs ($\mathbf{v_3}$,$\mathbf{v_{11}}$) and ($\mathbf{v_4}$,$\mathbf{v_{12}}$). One aspect of II-OEE flexible suspensions that is not apparent from the image is that the solid angles of the apical vertexes are constant and equal; $\sigma_u=\sigma_w=2\pi$.

**Example III-OAE, III-OAS Suspensions:** Parameter sets for the sub-types III-OAE and III-OAS are provided in Tables D-VI and D-V respectively. In addition to the parameter values the value of L (Eqs. 7) and K (Eqs. 9) and of the number of distinct Flexible Suspensions (FS) that have been created from the parameter set are included.

| # | N | L | Parameters | FS Count |
|---|---|---|---|---|
| **1** | **4** | **3** | **$l_1$=10, $\alpha_1$=45, $\beta_1$=55, $\alpha_2$=30, $\beta_2$=20.** | **3** |
| **2** | **6** | **3** | **$l_1$=10, $\alpha_1$=45, $\beta_1$=55, $\alpha_2$=30, $\beta_2$=20, $\alpha_3$=25.** | **3** |
| **3** | **6** | **4** | **$l_1$=10, $\alpha_1$=22, $\beta_1$=70, $\alpha_2$=100, $\beta_2$=45, $\alpha_3$=40.** | **1** |
| **4** | **8** | **5** | **$l_1$=10, $\alpha_1$=69, $\beta_1$=55, $\alpha_2$=30, $\beta_2$=45, $\alpha_3$=25, $\alpha_5$=21.** | **3** |
| **5** | **10** | **5** | **$l_1$=10, $\alpha_1$=54, $\beta_1$=75, $\alpha_2$=47, $\beta_2$=47, $\alpha_3$=35, $\alpha_5$=6, $\alpha_7$=41.** | **8** |
| **6** | **12** | **7** | **$l_1$=10, $\alpha_1$=75, $\beta_1$=70, $\alpha_2$=30, $\beta_2$=50, $\alpha_3$=25, $\alpha_5$=35, $\alpha_7$=20, $\alpha_9$=40.** | **5** |
| **7** | **14** | **8** | **$l_1$=10, $\alpha_1$=71, $\beta_1$=76, $\alpha_2$=70, $\beta_2$=31, $\alpha_3$=32, $\alpha_5$=29, $\alpha_7$=23, $\alpha_9$=37, $\alpha_{11}$=30.** | **13** |
| **8** | **16** | **8** | **$l_1$=10, $\alpha_1$=71, $\beta_1$=76, $\alpha_2$=70, $\beta_2$=30, $\alpha_3$=36, $\alpha_5$=33, $\alpha_7$=23, $\alpha_9$=37, $\alpha_{11}$=31, $\alpha_{13}$=28.** | **10** |

**Table D-IV.** III-OAE Suspension Example Parameter Sets

| # | N | K | Parameters | FS Count |
|---|---|---|---|---|
| **1** | **4** | **1** | **$l_1$=10, $\alpha_1$=105, $\beta_1$=30, $\alpha_2$=110, $\beta_2$=25.** | **2** |
| **2** | **6** | **1** | **$l_1$=10, $\alpha_1$=45, $\beta_1$=80, $\alpha_2$=42, $\beta_2$=37, $\alpha_3$=48.** | **4** |
| **3** | **8** | **1** | **$l_1$=10, $\alpha_1$=30, $\beta_1$=85, $\alpha_2$=40, $\beta_2$=70, $\alpha_3$=45, $\alpha_5$=45.** | **5** |
| **4** | **10** | **1** | **$l_1$=10, $\alpha_1$=30, $\beta_1$=91, $\alpha_2$=31, $\beta_2$=75, $\alpha_3$=27, $\alpha_5$=25, $\alpha_7$=30.** | **3** |
| **5** | **12** | **1** | **$l_1$=10, $\alpha_1$=41, $\beta_1$=104, $\alpha_2$=29, $\beta_2$=75, $\alpha_3$=27, $\alpha_5$=48, $\alpha_7$=21, $\alpha_9$=33.** | **1** |
| **6** | **14** | **2** | **$l_1$=10, $\alpha_1$=55, $\beta_1$=55, $\alpha_2$=100, $\beta_2$=40, $\alpha_3$=49, $\alpha_5$=57, $\alpha_7$=55, $\alpha_9$=48, $\alpha_{11}$=50.** | **6** |
| **7** | **16** | **3** | **$l_1$=10, $\alpha_1$=70, $\beta_1$=63, $\alpha_2$=69, $\beta_2$=59, $\alpha_3$=62, $\alpha_5$=66, $\alpha_7$=69, $\alpha_9$=67, $\alpha_{11}$=68, $\alpha_{13}$=65.** | **3** |

**Table D-V.** III-OAS Suspension Example Parameter Sets

The creation of valid complete parameter sets is a major computational challenge associated with these two sub-types. At each stage of computation (Appendix B) the computed angles may be invalid thereby invalidating the results of previous stages and necessitating an iteration on values that are satisfactory at previous stages. At this time we have performed this by a combination of analysis and experimentation and have not as yet attempted to automate this via software; however once a valid complete set has been established the determination of other complete sets from the M+3 parameter set is performed by an automated procedure.

An **III-OAE Flexible Suspension** is shown is Figure D-4. This suspension is constructed from the parameter set #8 in Table D-IV, has 32 faces (N=16) and has been created with the variable of flexion $\varepsilon_1$=75$^o$. The DI value for this is 3697=E71$_{16}$. Fig. D-4A is a front view of the entire polyhedra that shows the circular motion traces of the apical vertex **w** and the plane in which the fixed face $\mathbf{uv_2v_1}$ is located. The insets are of the same view without faces and show two meshes that are formed by two edges under flexion. These illustrate the symmetry that exists between motion with the variable of flexion $\varepsilon_1$ in the interval [0$^o$,180$^o$] and [180$^o$,360$^o$]. Fig. D-4B is top down view of the polyhedra framework that illustrates the flat co-planar folding defined by the second of equations (7); ie. $\delta_1=\delta_L=\pi$ ,

L=8, while all other dihedrals at the apical vertexes are =0. Insets show the specific edges that have this open folding.

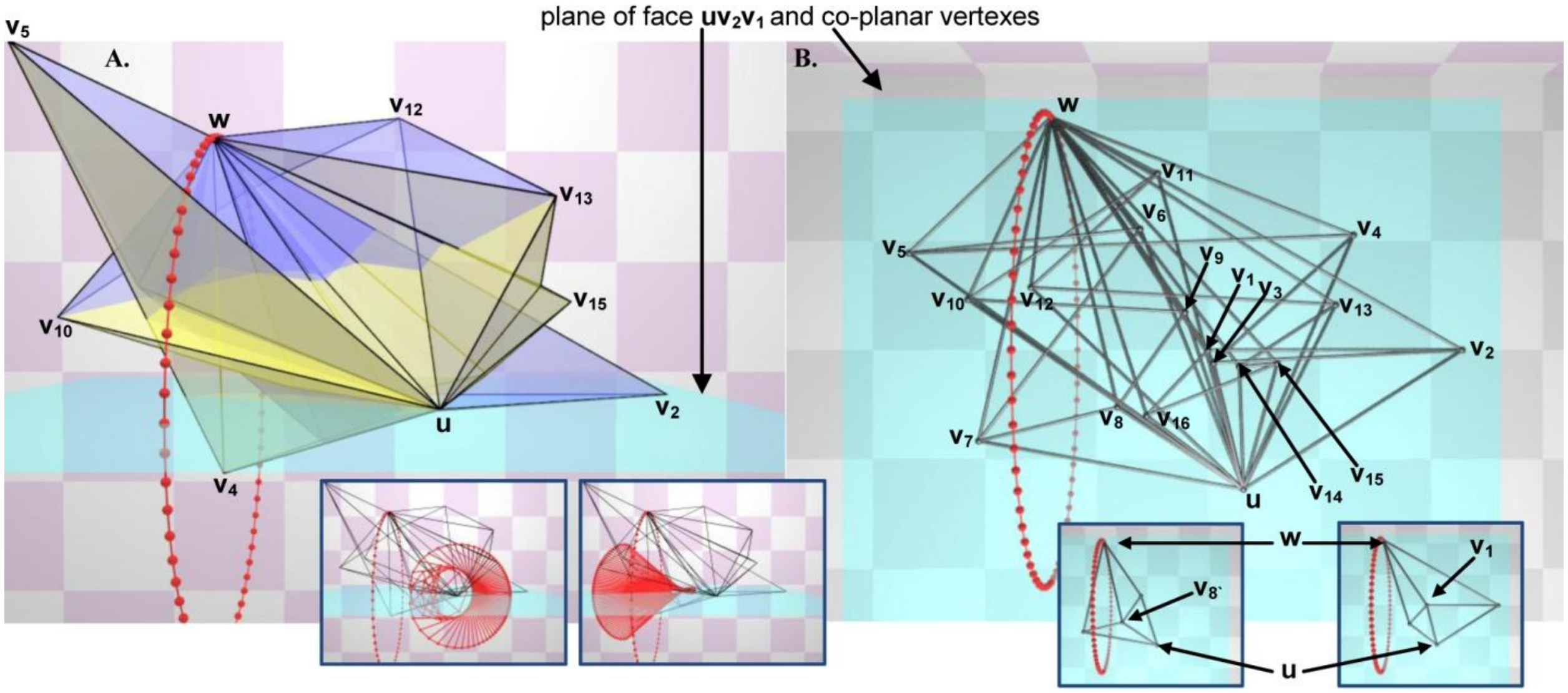

**Figure D-4.** III-OAE Flexible Suspension. **A.** Front View. **B.** Top View.

An **III-OAS Flexible Suspension** is illustrated is Figure D-5. This suspension is constructed from the parameter set #7 in Table D-V, has 32 faces (N=16) and a DI value of 14549=38D5$_{16}$.

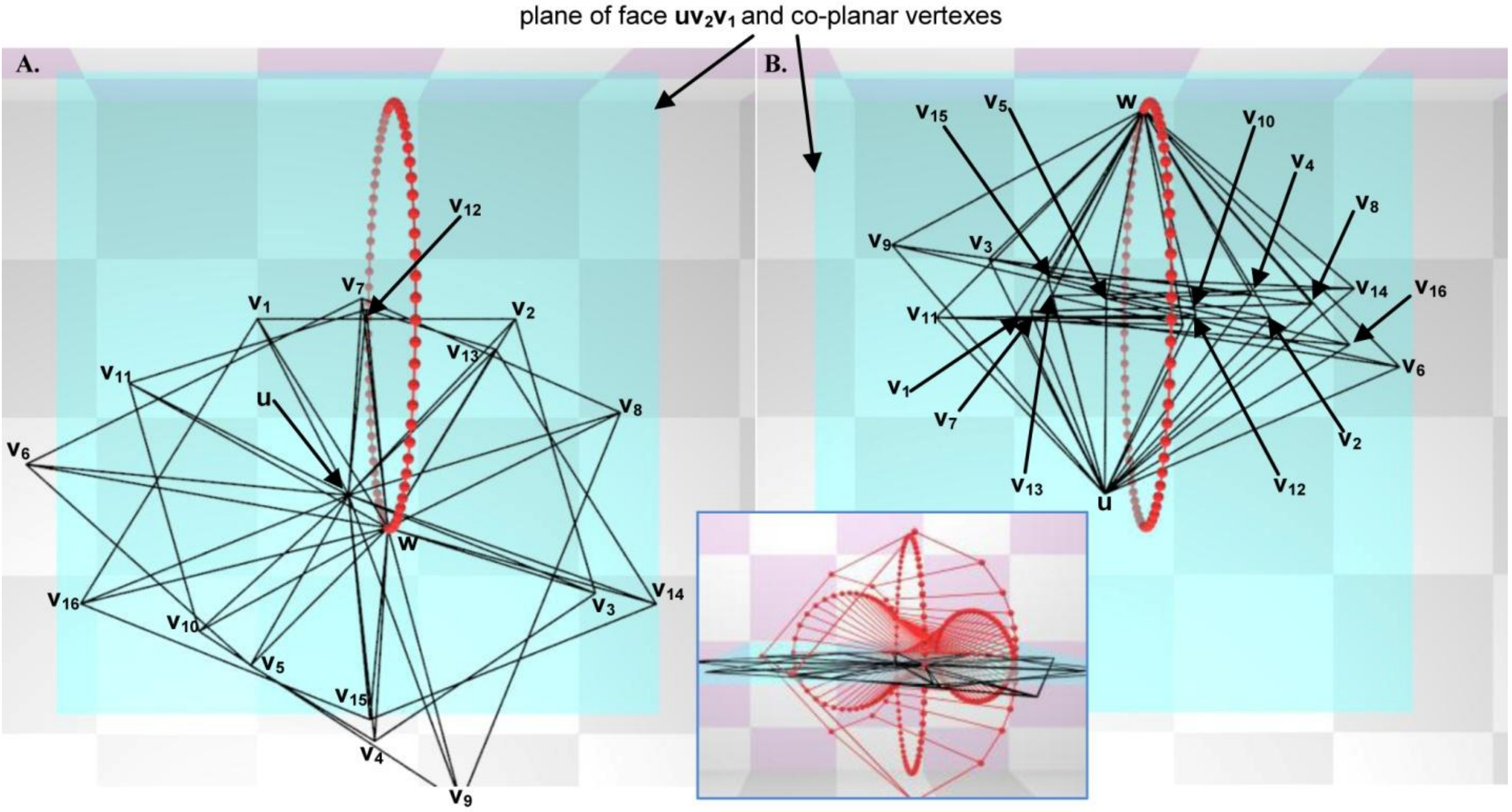

**Figure D-5.** III-OAS Flexible Suspension. **A.** Top View. **B.** Top View.

Fig. D-5A is a top down view of the polyhedra that shows co-planar framework of the edges with the variable of flexion $\varepsilon_1 = 0^o$ along with the circular trace of the motion of the apical vertex **w.** In this flat position all dihedral angles at the apical vertexes are $=\pi$. Because the parameter K=3, traversal of non-apical vertexes in sequential order involves three circuits encircling the apical vertexes. Fig. D-5B is a top down view that shows co-planar framework with the variable of flexion $\varepsilon_1 = 180^o$**.** In this flat position all dihedral angles at the apical vertexes are =0. The inset is a front view that shows the symmetrical mesh that is formed by two edges under flexion. These illustrate the symmetry that exists between motion with the variable of flexion $\varepsilon_1$ in the interval [$0^o$, $180^o$] and [$180^o$, $360^o$].

Computed parameters (face angles) for selected III-OAE and III-OAS examples are shown in Tables D-VI and D-VII respectively. Also shown is the DI associated with the computed parameter set. These tables when combined with the independent parameters in Tables D-IV and D-V define sufficient data to construct the example flexible polyhedra.

| # | N | Computed Angles | DI, $DI_{16}$ |
|---|---|---|---|
| **1** | 4 | **$\alpha_4$=30; $\beta_3$=112.12184, $\beta_4$=87.83527.** | **9,$9_{16}$** |
| **2** | 6 | **$\alpha_4$=16.40308, $\alpha_6$=13.59692; $\beta_3$=112.12184, $\beta_4$=118.35370, $\beta_5$=126.47742, $\beta_6$=104.23836.** | **41,$29_{16}$** |
| **3** | 6 | **$\alpha_4$=83.30367, $\alpha_6$=16.69633; $\beta_3$=21.45580, $\beta_4$=60.58002, $\beta_5$=74.47896, $\beta_6$=82.03499.** | **53,$35_{16}$** |
| **4** | 8 | **$\alpha_4$=83.32691, $\alpha_6$=60.25878, $\alpha_8$=53.06814;<br>$\beta_3$=107.57407, $\beta_4$=60.98961, $\beta_5$=111.23558, $\beta_6$=23.65624, $\beta_7$=79.35616, $\beta_8$=41.30176.** | 55,$37_{16}$ |
| **5** | 10 | **$\alpha_4$=100.88503, $\alpha_6$=94.89234, $\alpha_8$=42.70820, $\alpha_{10}$=10.28448;<br>$\beta_3$=27.48350, $\beta_4$=53.22171, $\beta_5$=151.69262, $\beta_6$= 22.74071, $\beta_7$=126.67290, $\beta_8$=37.49349, $\beta_9$=100.91108, $\beta_{10}$=12.47047.** | **805,$325_{16}$** |
| **6** | **12** | **$\alpha_4$=64.27757, $\alpha_6$=66.40444, $\alpha_8$=67.23134, $\alpha_{10}$=63.46828, $\alpha_{12}$=29.98240;<br>$\beta_3$=113.98589, $\beta_4$=42.89481, $\beta_5$=53.01637, $\beta_6$=69.15124, $\beta_7$=140.43696, $\beta_8$=25.11212, $\beta_9$=111.18488, $\beta_{10}$=36.54932, $\beta_{11}$=65.99163, $\beta_{12}$=26.93834.** | **3149, $C4D_{16}$** |
| **7** | 14 | **$\alpha_4$=79.94471, $\alpha_6$=64.21691, $\alpha_8$=28.13814, $\alpha_{10}$=69.71837, $\alpha_{12}$=106.32765, $\alpha_{14}$=9.97746 ;<br>$\beta_3$=130.35037, $\beta_4$=72.30741, $\beta_5$=29.55757, $\beta_6$=46.24987, $\beta_7$=18.56180, $\beta_8$=90.38385, $\beta_9$=131.26467, $\beta_{10}$=51.31861, $\beta_{11}$=24.06035, $\beta_{12}$=41.67901, $\beta_{13}$=19.54221, $\beta_{14}$=122.67419.** | **1009, $3F1_{16}$** |
| 8 | 16 | **$\alpha_4$=77.23450, $\alpha_6$=58.22912, $\alpha_8$=38.51987, $\alpha_{10}$=40.02534, $\alpha_{12}$=72.02307, $\alpha_{14}$=16.65480, $\alpha_{16}$=38.24054;<br>$\beta_3$=119.99899, $\beta_4$=64.61507, $\beta_5$=41.84098, $\beta_6$=57.99019, $\beta_7$=24.38078, $\beta_8$=88.80814, $\beta_9$=111.30147, $\beta_{10}$=69.20590, $\beta_{11}$=62.35174, $\beta_{12}$=40.78832, $\beta_{13}$=30.57510, $\beta_{14}$=93.17456, $\beta_{15}$=53.05563, $\beta_{16}$=85.23613.** | **3697, $E71_{16}$** |

**Table D-VI.** III-OAE Computed Angles

| # | N | Computed Angles | DI, $DI_{16}$ |
|---|---|---|---|
| **1** | 4 | **$\alpha_4$=70.00000; $\beta_3$=82.95205, $\beta_4$=50.47518.** | **3,$3_{16}$** |
| **2** | 6 | **$\alpha_4$=105.95429, $\alpha_6$=32.04571; $\beta_3$=99.13918, $\beta_4$=49.21957, $\beta_5$=37.80644, $\beta_6$=28.97738.** | **49,$31_{16}$** |
| **3** | 8 | **$\alpha_4$=91.18022, $\alpha_6$=47.37676, $\alpha_8$=1.44302;<br>$\beta_3$=40.43044, $\beta_4$=50.99966, $\beta_5$=34.62276, $\beta_6$=120.10265, $\beta_7$=27.97615, $\beta_8$=127.09938.** | **53,$35_{16}$** |
| **4** | 10 | **$\alpha_4$=2.29641, $\alpha_6$=65.48725, $\alpha_8$=56.38952, $\alpha_{10}$=24.82681;<br>$\beta_3$=105.25131, $\beta_4$=84.73441, $\beta_5$=130.85701, $\beta_6$=48.10408, $\beta_7$=63.33243, $\beta_8$=56.76065, $\beta_9$=35.09934, $\beta_{10}$=65.39581.** | **799,$31F_{16}$** |
| **5** | **12** | **$\alpha_4$=90.80897, $\alpha_6$=16.13559, $\alpha_8$=3.70555, $\alpha_{10}$=23.68880, $\alpha_{12}$=16.66108;<br>$\beta_3$=53.32442, $\beta_4$=35.56317, $\beta_5$=23.46576, $\beta_6$=104.33625, $\beta_7$=67.20918, $\beta_8$=115.40075, $\beta_9$=87.09237, $\beta_{10}$=109.81759, $\beta_{11}$=109.93166, $\beta_{12}$=129.25846.** | **4081,$FF1_{16}$** |
| **6** | 14 | **$\alpha_4$=52.66573, $\alpha_6$=17.47293, $\alpha_8$=9.54257, $\alpha_{10}$=25.82052, $\alpha_{12}$=137.95959, $\alpha_{14}$=16.53867;<br>$\beta_3$=34.84383, $\beta_4$=111.74521, $\beta_5$=29.50679, $\beta_6$=120.94430, $\beta_7$=53.40928, $\beta_8$=47.83117, $\beta_9$=91.51635, $\beta_{10}$=21.79280, $\beta_{11}$=108.12969, $\beta_{12}$=19.54775, $\beta_{13}$=39.29289, $\beta_{14}$=101.21476.** | **6197,$1835_{16}$** |
| 7 | 16 | **$\alpha_4$=35.52177, $\alpha_6$=88.31987, $\alpha_8$=87.36781, $\alpha_{10}$=70.20535, $\alpha_{12}$=39.37543, $\alpha_{14}$=64.82772, $\alpha_{16}$=85.38206;<br>$\beta_3$=56.25741, $\beta_4$=52.29144, $\beta_5$=77.21017, $\beta_6$=32.20233, $\beta_7$=69.65471, $\beta_8$=50.95207, $\beta_9$=36.56857, $\beta_{10}$=62.72476, $\beta_{11}$=44.85946, $\beta_{12}$=78.08803, $\beta_{13}$=76.63756, $\beta_{14}$=42.06629, $\beta_{15}$=62.48605, $\beta_{16}$=36.20026.** | **14549,$38D5_{16}$** |

**Table D-VII.** III-OAS Computed Angles

**Rigid Isomers:** Associated with each of the parameter sets for flexible suspensions are a number of rigid isomers that have the same faces and combinatorial structure. While we have not studied this topic in any detail it appears that a large number of isomers exist for any given parameter set or flexible suspension. For each flexible folding of a sub-type I-OEE suspension there are another $2^{M-1}$ rigid isomers. Each has a specific $\varepsilon_1$ and an DI identifier value that is appropriate for its construction. Similarly for sub-types II-AEE and II-OEE. On the other hand for sub-types III-OAE and III-OAS, all foldings exhibit isomers at $\varepsilon_1=0$ or $\varepsilon_1=\pi$ which are identical to the flexible folding evaluated at $\varepsilon_1=0$ or $\varepsilon_1=\pi$.